\font\tenmib=cmmib10
\font\sevenmib=cmmib10 scaled 800
\font\titolo=cmbx12
\font\titolone=cmbx10 scaled\magstep 2

\font\cs=cmcsc10
\font\sc=cmcsc10
\font\css=cmcsc8

\font\ninerm=cmr9
\font\ottorm=cmr8
\textfont5=\tenmib\scriptfont5=\sevenmib\scriptscriptfont5=\fivei

\font\euftw=eufm10 scaled\magstep1
\font\msytw=msbm10 scaled\magstep1

\font\msytwww=msbm7 scaled\magstep1
\font\indbf=cmbx10 scaled\magstep2

\font\ottorm=cmr8\font\ottoi=cmmi8\font\ottosy=cmsy8%
\font\ottobf=cmbx8\font\ottott=cmtt8
\font\ottocss=cmcsc8%
\font\ottosl=cmsl8\font\ottoit=cmti8%
\font\sixrm=cmr6\font\sixbf=cmbx6\font\sixi=cmmi6\font\sixsy=cmsy6%
\font\fiverm=cmr5\font\fivesy=cmsy5\font\fivei=cmmi5\font\fivebf=cmbx5%

\def\ottopunti{\def\rm{\fam0\ottorm}%
\textfont0=\ottorm\scriptfont0=\sixrm\scriptscriptfont0=\fiverm%
\textfont1=\ottoi\scriptfont1=\sixi\scriptscriptfont1=\fivei%
\textfont2=\ottosy\scriptfont2=\sixsy\scriptscriptfont2=\fivesy%
\textfont3=\tenex\scriptfont3=\tenex\scriptscriptfont3=\tenex%
\textfont4=\ottocss\scriptfont4=\sc\scriptscriptfont4=\sc%
\textfont\itfam=\ottoit\def\it{\fam\itfam\ottoit}%
\textfont\slfam=\ottosl\def\sl{\fam\slfam\ottosl}%
\textfont\ttfam=\ottott\def\tt{\fam\ttfam\ottott}%
\textfont\bffam=\ottobf\scriptfont\bffam=\sixbf%
\scriptscriptfont\bffam=\fivebf\def\bf{\fam\bffam\ottobf}%
\setbox\strutbox=\hbox{\vrule height7pt depth2pt width0pt}%
\normalbaselineskip=9pt\let\sc=\sixrm\normalbaselines\rm}
\let\nota=\ottopunti%

%
%
%
%
%
%
%

\global\newcount\numsec\global\newcount\numapp
\global\newcount\numfor\global\newcount\numfig
\global\newcount\numsub
\numsec=0\numapp=0\numfig=0
\def\veroparagrafo{\number\numsec}\def\veraformula{\number\numfor}
\def\veraappendice{\number\numapp}\def\verasub{\number\numsub}
\def\verafigura{\number\numfig}

\def\section(#1,#2){\advance\numsec by 1\numfor=1\numsub=1\numfig=1%
\SIA p,#1,{\veroparagrafo} %
\write15{\string\Fp (#1){\secc(#1)}}%
\write16{ sec. #1 ==> \secc(#1)  }%
\hbox to \hsize{\titolo\hfill \number\numsec. #2\hfill%
\expandafter{\alato(sec. #1)}}\*}

\def\appendix(#1,#2){\advance\numapp by 1\numfor=1\numsub=1\numfig=1%
\SIA p,#1,{A\veraappendice} %
\write15{\string\Fp (#1){\secc(#1)}}%
\write16{ app. #1 ==> \secc(#1)  }%
\hbox to \hsize{\titolo\hfill Appendix A\number\numapp. #2\hfill%
\expandafter{\alato(app. #1)}}\*}

\def\senondefinito#1{\expandafter\ifx\csname#1\endcsname\relax}

\def\SIA #1,#2,#3 {\senondefinito{#1#2}%
\expandafter\xdef\csname #1#2\endcsname{#3}\else
\write16{???? ma #1#2 e' gia' stato definito !!!!} \fi}

\def \Fe(#1)#2{\SIA fe,#1,#2 }
\def \Fp(#1)#2{\SIA fp,#1,#2 }
\def \Fg(#1)#2{\SIA fg,#1,#2 }

\def\etichetta(#1){(\veroparagrafo.\veraformula)%
\SIA e,#1,(\veroparagrafo.\veraformula) %
\global\advance\numfor by 1%
\write15{\string\Fe (#1){\equ(#1)}}%
\write16{ EQ #1 ==> \equ(#1)  }}

\def\etichettaa(#1){(A\veraappendice.\veraformula)%
\SIA e,#1,(A\veraappendice.\veraformula) %
\global\advance\numfor by 1%
\write15{\string\Fe (#1){\equ(#1)}}%
\write16{ EQ #1 ==> \equ(#1) }}

\def\getichetta(#1){\veroparagrafo.\verafigura%
\SIA g,#1,{\veroparagrafo.\verafigura} %
\global\advance\numfig by 1%
\write15{\string\Fg (#1){\graf(#1)}}%
\write16{ Fig. #1 ==> \graf(#1) }}

\def\etichettap(#1){\veroparagrafo.\verasub%
\SIA p,#1,{\veroparagrafo.\verasub} %
\global\advance\numsub by 1%
\write15{\string\Fp (#1){\secc(#1)}}%
\write16{ par #1 ==> \secc(#1)  }}

\def\etichettapa(#1){A\veraappendice.\verasub%
\SIA p,#1,{A\veraappendice.\verasub} %
\global\advance\numsub by 1%
\write15{\string\Fp (#1){\secc(#1)}}%
\write16{ par #1 ==> \secc(#1)  }}

\def\Eq(#1){\eqno{\etichetta(#1)\alato(#1)}}
\def\eq(#1){\etichetta(#1)\alato(#1)}
\def\Eqa(#1){\eqno{\etichettaa(#1)\alato(#1)}}
\def\eqa(#1){\etichettaa(#1)\alato(#1)}
\def\eqg(#1){\getichetta(#1)\alato(fig. #1)}
\def\sub(#1){\0\palato(p. #1){\bf \etichettap(#1).}}
\def\asub(#1){\0\palato(p. #1){\bf \etichettapa(#1).}}

\def\equv(#1){\senondefinito{fe#1}$\clubsuit$#1%
\write16{eq. #1 non e' (ancora) definita}%
\else\csname fe#1\endcsname\fi}
\def\grafv(#1){\senondefinito{fg#1}$\clubsuit$#1%
\write16{fig. #1 non e' (ancora) definito}%
\else\csname fg#1\endcsname\fi}
\def\secv(#1){\senondefinito{fp#1}$\clubsuit$#1%
\write16{par. #1 non e' (ancora) definito}%
\else\csname fp#1\endcsname\fi}

\def\equ(#1){\senondefinito{e#1}\equv(#1)\else\csname e#1\endcsname\fi}
\def\graf(#1){\senondefinito{g#1}\grafv(#1)\else\csname g#1\endcsname\fi}
\def\figura(#1){{\css Figure} \getichetta(#1)}
\def\secc(#1){\senondefinito{p#1}\secv(#1)\else\csname p#1\endcsname\fi}
\def\sec(#1){{\S\secc(#1)}}
\def\refe(#1){{[\secc(#1)]}}

\def\BOZZA{\boz=1
\def\alato(##1){\rlap{\kern-\hsize\kern-1.2truecm{$\scriptstyle##1$}}}
\def\palato(##1){\rlap{\kern-1.2truecm{$\scriptstyle##1$}}}
}

\def\alato(#1){}
\def\galato(#1){}
\def\palato(#1){}


{\count255=\time\divide\count255 by 60 \xdef\hourmin{\number\count255}
        \multiply\count255 by-60\advance\count255 by\time
   \xdef\hourmin{\hourmin:\ifnum\count255<10 0\fi\the\count255}}

\def\oramin{\hourmin }

\def\data{\number\day/\ifcase\month\or gennaio \or febbraio \or marzo \or
aprile \or maggio \or giugno \or luglio \or agosto \or settembre
\or ottobre \or novembre \or dicembre \fi/\number\year;\ \oramin}
\footline={\rlap{\hbox{\copy200}}\tenrm\hss \number\pageno\hss}


\newcount\driver 
\newdimen\xshift \newdimen\xwidth
\def\ins#1#2#3{\vbox to0pt{\kern-#2 \hbox{\kern#1
#3}\vss}\nointerlineskip}

\def\insertplot#1#2#3#4#5{\par%
\xwidth=#1 \xshift=\hsize \advance\xshift by-\xwidth \divide\xshift by 2%
\yshift=#2 \divide\yshift by 2%
\line{\hskip\xshift \vbox to #2{\vfil%
\ifnum\driver=0 #3
\special{ps: plotfile #4.eps} 
\fi
\ifnum\driver=1 #3 \includegraphics{#4.eps}\fi
\ifnum\driver=2 #3
\ifnum\mgnf=0\special{#4.eps 1. 1. scale} \fi
\ifnum\mgnf=1\special{#4.eps 1.2 1.2 scale}\fi
\fi }\hfill \raise\yshift\hbox{#5}}}

\def\insertplotttt#1#2#3#4{\par%
\xwidth=#1 \xshift=\hsize \advance\xshift by-\xwidth \divide\xshift by 2%
\yshift=#2 \divide\yshift by 2%
\line{\hskip\xshift \vbox to #2{\vfil%
\ifnum\driver=0 #3
\special{ps: plotfile #4.eps} 
\fi
\ifnum\driver=1 #3 \includegraphics{#4.eps}\fi
\ifnum\driver=2 #3
\ifnum\mgnf=0\special{#4.eps 1. 1. scale} \fi
\ifnum\mgnf=1\special{#4.eps 1.2 1.2 scale}\fi
\fi }\hfill}}

\def\insertplott#1#2#3{\par%
\xwidth=#1 \xshift=\hsize \advance\xshift by-\xwidth \divide\xshift by 2%
\yshift=#2 \divide\yshift by 2%
\line{\hskip\xshift \vbox to #2{\vfil%
\ifnum\driver=0
\special{ps: plotfile #3.eps} 
 \fi
\ifnum\driver=1 \includegraphics{#3.eps}\fi
\ifnum\driver=2
\ifnum\mgnf=0\special{#3.eps 1. 1. scale} \fi 
\ifnum\mgnf=1\special{#3.eps 1.2 1.2 scale}\fi
\fi }\hfill}}

\newdimen\xshift \newdimen\xwidth \newdimen\yshift
\def\eqfig#1#2#3#4#5{
\par\xwidth=#1 \xshift=\hsize \advance\xshift
by-\xwidth \divide\xshift by 2
\yshift=#2 \divide\yshift by 2
\line{\hglue\xshift \vbox to #2{\vfil
\ifnum\driver=0 #3
\special{ps: plotfile #4.ps} 
\fi
\ifnum\driver=1 #3 \includegraphics{#4.ps}\fi
\ifnum\driver=2 #3 \special{
\ifnum\mgnf=0 #4.ps 1. 1. scale \fi
\ifnum\mgnf=1 #4.ps 1.2 1.2 scale\fi}
\fi}\hfill\raise\yshift\hbox{#5}}}

\let\a=\alpha \let\b=\beta  \let\g=\gamma  \let\d=\delta \let\e=\varepsilon
     \let\th=\theta  \let\l=\lambda
\let\m=\mu    \let\n=\nu    \let\x=\xi     \let\p=\pi    
\let\s=\sigma \let\t=\tau   \let\f=\varphi 
   \let\o=\omega

\def\\{\hfill\break} \let\==\equiv

\let\io=\infty 

\let\0=\noindent

\let\dpr=\partial

\def\tende#1{\,\vtop{\ialign{##\crcr\rightarrowfill\crcr
 \noalign{\kern-1pt\nointerlineskip}
 \hskip3.pt${\scriptstyle #1}$\hskip3.pt\crcr}}\,}
\def\otto{\,{\kern-1.truept\leftarrow\kern-5.truept\to\kern-1.truept}\,}

\def\MM{{\cal M}} \def\VV{{\cal V}}
\def\CC{{\cal C}}
\def\TT{{\cal T}}\def\BB{{\cal B}}
\def\RR{{\cal R}} 
\def\DD{{\cal D}}\def\SS{{\cal S}}

\def\T#1{{#1_{\kern-3pt\lower7pt\hbox{$\widetilde{}$}}\kern3pt}}
\def\VVV#1{{\underline #1}_{\kern-3pt
\lower7pt\hbox{$\widetilde{}$}}\kern3pt\,}
\def\W#1{#1_{\kern-3pt\lower7.5pt\hbox{$\widetilde{}$}}\kern2pt\,}

\def\lis{\overline}

\def\indica{\leaders \hbox to 0.5cm{\hss.\hss}\hfill}
\def\guida{\leaders\hbox to 1em{\hss.\hss}\hfill}

\def\V0{{\bf 0}} \def\BBBB{\hbox{\euftw B}}
\def\RRRR{\hbox{\euftw R}}


\mathchardef\aa   = "050B
\mathchardef\bb   = "050C
\mathchardef\xxx  = "0518
\mathchardef\hhh  = "0511
\mathchardef\zzzzz= "0510
\mathchardef\oo   = "0521
\mathchardef\lll  = "0515
\mathchardef\mm   = "0516
\mathchardef\Dp   = "0540
\mathchardef\H    = "0548
\mathchardef\FFF  = "0546
\mathchardef\ppp  = "0570
\mathchardef\nn   = "0517
\mathchardef\pps  = "0520
\mathchardef\FFF  = "0508
\mathchardef\nnnnn= "056E

\def\to{\rightarrow}

\let\ciao=\bye
\def\qed{\hfill\raise1pt\hbox{\vrule height5pt width5pt depth0pt}}

 \def\Val{{\rm Val}}
\def\indic{\hbox{\raise-2pt \hbox{\indbf 1}}}

\def\RRR{\hbox{\msytw R}} 
 \def\CCC{\hbox{\msytw C}}
 
\def\NNN{\hbox{\msytw N}} 
 \def\ZZZ{\hbox{\msytw Z}}
 \def\zzz{\hbox{\msytwww Z}}


\newcount\mgnf  
\mgnf=0 

\ifnum\mgnf=0
\def\openone{\leavevmode\hbox{\ninerm 1\kern-3.3pt\tenrm1}}%
\def\*{\vglue0.3truecm}\fi
\ifnum\mgnf=1
\def\openone{\leavevmode\hbox{\ninerm 1\kern-3.63pt\tenrm1}}%
\def\*{\vglue0.5truecm}\fi


\newcount\tipobib\newcount\boz\boz=0\newcount\aux\aux=1
\newdimen\bibskip\newdimen\maxit\maxit=0pt


\tipobib=0
\def\9#1{\ifnum\aux=1#1\else\relax\fi}

\newwrite\bib
\immediate\openout\bib=\jobname.bib
\global\newcount\bibex
\bibex=0
\def\verabib{\number\bibex}

\ifnum\tipobib=0
\def\cita#1{\expandafter\ifx\csname c#1\endcsname\relax
\hbox{$\clubsuit$}#1\write16{Manca #1 !}%
\else\csname c#1\endcsname\fi}
\def\rife#1#2#3{\immediate\write\bib{\string\raf{#2}{#3}{#1}}
\immediate\write15{\string\C(#1){[#2]}}
\setbox199=\hbox{#2}\ifnum\maxit < \wd199 \maxit=\wd199\fi}
\else
\def\cita#1{%
\expandafter\ifx\csname d#1\endcsname\relax%
\expandafter\ifx\csname c#1\endcsname\relax%
\hbox{$\clubsuit$}#1\write16{Manca #1 !}%
\else\probib(ref. numero )(#1)%
\csname c#1\endcsname%
\fi\else\csname d#1\endcsname\fi}%
\def\rife#1#2#3{\immediate\write15{\string\Cp(#1){%
\string\immediate\string\write\string\bib{\string\string\string\raf%
{\string\verabib}{#3}{#1}}%
\string\Cn(#1){[\string\verabib]}%
\string\CCc(#1)%
}}}%
\fi

\def\Cn(#1)#2{\expandafter\xdef\csname d#1\endcsname{#2}}
\def\CCc(#1){\csname d#1\endcsname}
\def\probib(#1)(#2){\global\advance\bibex+1%
\9{\immediate\write16{#1\verabib => #2}}%
}

\def\C(#1)#2{\SIA c,#1,{#2}}
\def\Cp(#1)#2{\SIAnx c,#1,{#2}}

\def\SIAnx #1,#2,#3 {\senondefinito{#1#2}%
\expandafter\def\csname#1#2\endcsname{#3}\else%
\write16{???? ma #1,#2 e' gia' stato definito !!!!}\fi}

\bibskip=10truept
\def\hboxto{\hbox to}

\catcode`\{=12\catcode`\}=12
\catcode`\<=1\catcode`\>=2
\immediate\write\bib<
        \string\halign{\string\hboxto \string\maxit%
        {\string #\string\hfill}&%
        \string\vtop{\string\parindent=0pt\string\advance\string\hsize%
        by -.5truecm%
        \string#\string\vskip \bibskip
        }\string\cr%
>
\catcode`\{=1\catcode`\}=2
\catcode`\<=12\catcode`\>=12

\def\raf#1#2#3{\ifnum \boz=0 [#1]&#2 \cr\else
\llap{${}_{\rm #3}$}[#1]&#2\cr\fi}

\newread\bibin

\catcode`\{=12\catcode`\}=12
\catcode`\<=1\catcode`\>=2
\def\chiudibib<
\catcode`\{=12\catcode`\}=12
\catcode`\<=1\catcode`\>=2
\immediate\write\bib<}>
\catcode`\{=1\catcode`\}=2
\catcode`\<=12\catcode`\>=12
>
\catcode`\{=1\catcode`\}=2
\catcode`\<=12\catcode`\>=12

\def\makebiblio{
\ifnum\tipobib=0
\advance \maxit by 10pt
\else
\maxit=1.truecm
\fi
\chiudibib
\immediate \closeout\bib
\openin\bibin=\jobname.bib
\ifeof\bibin\relax\else
\raggedbottom
\input \jobname.bib
\fi
}

\openin13=#1.aux \ifeof13 \relax \else
\input #1.aux \closein13\fi
\openin14=\jobname.aux \ifeof14 \relax \else
\input \jobname.aux \closein14 \fi
\immediate\openout15=\jobname.aux

\def\biblio{\*\*\centerline{\titolo References}\*\nobreak\makebiblio}


\ifnum\mgnf=0
   \magnification=\magstep0
   \hsize=15.2truecm\vsize=20.8truecm\voffset1.truecm\hoffset.5truecm
   \parindent=0.3cm\baselineskip=0.45cm\fi
\ifnum\mgnf=1
   \magnification=\magstep1\hoffset=0.truecm
   \hsize=15.2truecm\vsize=20.8truecm
   \baselineskip=18truept plus0.1pt minus0.1pt \parindent=0.9truecm
   \lineskip=0.5truecm\lineskiplimit=0.1pt      \parskip=0.1pt plus1pt\fi


\mgnf=0   
\driver=1 
\openin14=\jobname.aux \ifeof14 \relax \else
\input \jobname.aux \closein14 \fi
\openout15=\jobname.aux



\centerline{\titolone Summation of divergent series and Borel summability}
\vskip.3truecm
\centerline{\titolone for strongly dissipative equations}
\vskip.3truecm
\centerline{\titolone with periodic or quasi-periodic forcing terms}

\vskip1.truecm
\centerline{{\titolo
G. Gentile${^{\dagger}}$,
M.V. Bartuccelli${^{\star}}$,
J.H.B. Deane${^{\star}}$}}
\vskip0.4truecm
\centerline{{}$^{\dagger}$ Dipartimento di Matematica,
Universit\`a di Roma Tre, Roma, I-00146, Italy}
\vskip.1truecm
\centerline{$^{\star}$ Department of Mathematics and Statistics,
University of Surrey, Guildford, GU2 7XH, UK}
\vskip1.truecm

\line{\vtop{
\line{\hskip1.1truecm\vbox{\advance \hsize by -2.3 truecm
\0{\cs Abstract.}
{\it We consider a class of second order ordinary differential equations
describing one-dimensional systems with a quasi-periodic analytic
forcing term and in the presence of damping.
As a physical application one can think of a
resistor-inductor-varactor circuit with a periodic
(or quasi-periodic) forcing function, even if the range of
applicability of the theory is much wider.
In the limit of large damping we look for quasi-periodic solutions
which have the same frequency vector of the forcing term, and we study
their analyticity properties in the inverse of the damping coefficient.
We find that already the case of periodic forcing terms is non-trivial,
as the solution is not analytic in a neighbourhood of the origin:
it turns out to be Borel-summable. In the case of
quasi-periodic forcing terms we need Renormalization Group
techniques in order to control the small divisors arising in the
perturbation series. We show the existence of a summation criterion
of the series in this case also, but, however, this 
can not be interpreted as Borel summability.
}} \hfill} }}

\vskip1.truecm
\section(1,Introduction)

\0Consider the ordinary differential equation
$$ \e \ddot x + \dot x + \e x^{2} = \e f(\oo t) ,
\Eq(1.1) $$
where $\oo\in\RRR^{d}$ is the frequency vector,
$f(\pps)$ is an analytic function,
$$ f(\pps) = \sum_{\nn\in \zzz^{d}} {\rm e}^{i\nn \cdot \pps} f_{\nn} ,
\Eq(1.2) $$
with average $\a=a^{2}$, with $a>0$ (hence $\langle f \rangle \=
f_{0}=\a$), and $\e>0$ is a real parameter. Here and henceforth
we denote with $\cdot$ the scalar product in $\RRR^{d}$.
By the analyticity assumption on $f$ there are
two strictly positive constants $F$ and $\x$ such that
one has $|f_{\nn}|\le F \, {\rm e}^{-\x|\nn|}$ for all $\nn\in\ZZZ^{d}$.

By writing $\g=1/\e$ the equation becomes
$$ \ddot x + \g \dot x + x^{2} = f(\oo t) , 
\Eq(1.3) $$
which describes a nonlinear electronic circuit,
known as resistor-inductor-varactor circuit,
subject to a quasi-periodic forcing function.
Taking $d=1$ and $f(\o t)=\a+\b \sin t$,
this equation has been studied in Ref.~\cita{DMGB}, where,
among other things, it has been found numerically that for $\g$ large
enough there exists only one attracting periodic orbit and the
corresponding period is $2\p/\o=2\p$, the same of the forcing term.
Furthermore one can prove analytically that such a
periodic orbit is the only one in a neighbourhood
of radius $O(1/\g)$ around the point $(a,0)$.

Here we give some further analytical support to such numerical findings.
In particular we show that, if we take as forcing term an
analytic periodic function,
$$ f(\psi) = \sum_{\n\in \zzz} {\rm e}^{i\n \psi} f_{\n} ,
\qquad f_{0} = \a > 0 ,
\Eq(1.4) $$
then for $\e$ small enough there is
a $2\p/\o$-periodic solution, but this is not analytic in $\e=1/\g$
in a neighbourhood the origin in the complex $\e$-plane.
We find that such a solution is Borel-summable.

We also show that by considering quasi-periodic forcing terms,
as in \equ(1.3), we still have a quasi-periodic solution
with the same frequency vector $\oo$ of the forcing term,
but we can only say in general that such a solution is analytic
in a domain with boundary crossing the origin.

Finally we shall see that considering more general nonlinearities
introduces no further difficulties, and equations like
$$ \ddot x + \g \dot x + g(x) = f(\oo t) , \qquad
\lim_{|x|\to\io} {|g(x)|\over |x|} = \io ,
\Eq(1.5) $$
with $g$ and $f$ both analytic in their arguments, can be dealt
with essentially in the same way. Simply, we have to impose
a non-degeneracy condition on the function $g$, which reads as
$$ \exists x_{0} \hbox{ such that } g(x_{0}) = f_{0} \hbox{ and }
g'(x_{0}) \neq 0 .
\Eq(1.6) $$
In the particular case of homogeneous $g(x)$,
that is $g(x)=\s x^{p}$, with $p\ge 2$ an integer and $\s\in\RRR$,
the condition is automatically satisfied if $p$ is odd
(for any value of $\s$), while it requires $\s f_{0}>0$
for $p$ even, -- as assumed in \equ(1.1).

The paper is organized as follows. For expository clearness
we start with the case of periodic forcing terms.
In Sections \secc(2) and \secc(3) we show that a periodic solution
with frequency $\o$ in the form of a formal power series in $\e$
(perturbation series) is well defined to all orders,
and it admits a natural graphical representation.
In Section \secc(4) we study further such a series,
and we see that there is strong evidence to expect that it diverges.
The best bounds that we are able to provide for the coefficients
grow as factorials. To obtain bounds which allow summability
of the perturbation series we have to perform a suitable summation
in order to give the series a meaning. This is done in Section \secc(5),
and the resummed series is found to represent a $2\p/\o$-periodic
solution which is Borel-summable in $\e$. To prove the latter property
we rely on Nevanlinna's improvement of Watson's theorem \cita{S}.
In Section \secc(6) we consider the case of quasi-periodic
forcing terms. We find that the perturbation series is well defined
if the frequency vector of the forcing term satisfies a
Diophantine condition, and, by using Renormalization Group techniques
in order to deal with the small divisors problem,
we find that the resummed series still converges to a
quasi-periodic solution, and it defines a function
analytic in a domain containing the origin in its boundary.
We shall see that the bounds we find do not allow us
any more to obtain Borel summability, unlike the case
of periodic forcing terms.
In Section \secc(7) we discuss how to extend the analysis to more
general nonlinearities $g(x)$, by requiring the condition \equ(1.6)
to be satisfied.

The interest of the approach we propose is that it allows
using perturbation theory which can be very natural in problems
in which a small parameter appears. In fact analyticity
in $\e$ for $\e$ close to $0$ (that is in $\g$ for $\g$ large enough)
could be proved very likely with other techniques, but a naive
expansion in powers of $\e$ is prevented by the lack of analyticity
in a neighbourhood of the origin.
On the other hand the perturbation series gives a very accurate
description of the solution, hence it is important to know that such a
series is an asymptotic series, and its use is fully justified.
Finally we can mention that the quasi-periodic solution we
investigate is of physical relevance, hence it can be useful
to study its properties. For instance in the case of the aforementioned
resistor-inductor-varactor circuit in Ref.~\cita{DMGB} the
$2\pi/\o$-periodic solution is numerically found to attract any
trajectory which remains bounded in phase space.

\*\*
\section(2,Formal analysis)

\0Consider first \equ(1.1) for $d=1$, that is
$$ \e \ddot x + \dot x + \e x^{2} = \e f(\o t) ,
\Eq(2.1) $$
with $f(\psi)$ given by \equ(1.4).
We look for bounded solutions (if any) which are
analytic in $\e$, that is of the form
$$ x(t) = \sum_{k=0}^{\io} \e^{k} x^{(k)}(t) .
\Eq(2.2) $$
Inserting \equ(2.2) into \equ(2.1) and equating terms of the same Taylor 
order we find the set of recursive equations
$$ \eqalign{
\dot x^{(0)} & = 0 , \cr
\dot x^{(1)} & = - \ddot x^{(0)} - x^{(0)2} + f , \cr
\dot x^{(k)} & = - \ddot x^{(k-1)} - \sum_{k_{1}+k_{2}=k-1}
x^{(k_{1})} x^{(k_{2})} , \qquad k \ge 2 . \cr}
\Eq(2.3) $$
From the first equation (zeroth order) we obtain that $x^{(0)}$ has to be
constant, say $x^{(0)}=c_{0}$ with $c_{0}$ to be determined.
The second equation (first order) can give a bounded solution
only if $-c_{0}^{2}+\a=0$, which fixes $c_{0}=\sqrt{\a}=a$
and gives $x^{(1)}(t)$ as a periodic function with the
same period of the forcing term:
$$ x^{(1)}(t) = x^{(1)}(0) +
\int_{0}^{t} {\rm d} t' \left( f(\o t') - \a \right) .
\Eq(2.4) $$
As each $x^{(k)}(t)$ depends on the functions $x^{(k')}(t)$
with $k'<k$, we expect that if there is any
periodic solution then it must have the same period as the forcing term.

To continue the analysis to all orders it is more convenient
to write the recursive equations \equ(2.3) in Fourier space.
The analysis to first order and the considerations above motivate
us to write in \equ(2.2)
$$ x(t) = \sum_{k=0}^{\io} \e^{k} x^{(k)}(t) = 
\sum_{k=0}^{\io} \e^{k} \sum_{\n\in\zzz} e^{i\n\o t} x^{(k)}_{\n} ,
\Eq(2.5) $$
which inserted into \equ(2.3) gives for $\n \neq 0$
$$ \eqalign{
x^{(0)}_{\n} & = 0 , \cr
x^{(1)}_{\n} & = {f_{\n} \over i\o\n} , \cr
x^{(k)}_{\n} & = - (i\o\n)\, x^{(k-1)}_{\n} - {1 \over i\o\n}
\sum_{\scriptstyle k_{1}+k_{2}=k-1 \atop k_{1},k_{2} \ge 0}
\sum_{\n_{1}+\n_{2}=\n}
x^{(k_{1})}_{\n_{1}} x^{(k_{2})}_{\n_{2}} , \qquad k \ge 2 , \cr}
\Eq(2.6) $$
provided that one has for $\n=0$
$$ \eqalign{
0 & = - x^{(0)2}_{0} + f_{0} , \cr
0 & = \sum_{\scriptstyle k_{1}+k_{2}=k \atop k_{1},k_{2} \ge 0}
\sum_{\n_{1}+\n_{2}=0}
x^{(k_{1})}_{\n_{1}} x^{(k_{2})}_{\n_{2}} , \qquad k \ge 1 . \cr}
\Eq(2.7) $$
If we set $x^{(k)}_{0}=c_{k}$ then the first of \equ(2.7) fixes,
as already noted,
$$ c_{0}=a=\sqrt{\a} ,
\Eq(2.8) $$
because one has $f_{0}=\a>0$, while the second of \equ(2.7) gives
$$ \sum_{k'=0}^{k} \sum_{\n_{1}\in\zzz} x^{(k-k')}_{\n_{1}}
x^{(k')}_{-\n_{1}} = 0 .
\Eq(2.9) $$
The latter equation, by taking into account \equ(2.8) and
the first of \equ(2.6), can be more conveniently written as
$$ c_{1} = 0 , \qquad\qquad
c_{k} = -{1 \over 2 c_{0}} \sum_{k'=1}^{k-1}
\sum_{\n_{1}\in\zzz} x^{(k-k')}_{\n_{1}} x^{(k')}_{-\n_{1}} ,
\qquad k \ge 2 ,
\Eq(2.10) $$
which provides an iterative definition of the coefficients
$c_{k}$ as the right hand side depends only on the coefficients
$c_{k'}$ with $k'<k$. To deduce $c_{1}=0$ we used the first of \equ(2.6),
which, inserted into \equ(2.9) for $k=1$, gives $2 c_{0}c_{1} = 0$,
hence $c_{1}=0$ as $c_{0} \neq 0$.

The following result holds.

\*

\0{\bf Lemma 1.} {\it Consider \equ(2.1) with $f$ given by \equ(1.4)
Then there exists a formal power series solution \equ(2.2)
whose coefficients $x^{(k)}(t)$ are analytic in $t$.
If $f$ is a trigonometric polynomial, that is
in \equ(1.4) one has $|\n|\le N$ for some $N\in\NNN$, then for all
$k \ge 0$ the functions $x^{(k)}(t)$ are trigonometric polynomials of
order $[(k+1)/2]N$, where $[\cdot]$ denotes the integer part. This means
that one has $x^{(2k)}_{\n}=0$ and $x^{(2k-1)}_{\n}=0$ for $|\n|>kN$.}

\*

\0{\it Proof.} The existence of a formal solution \equ(2.2),
with coefficients $x^{(k)}(t)$ analytic in $t$ for all $k\ge 0$,
follows from the analysis above. If $f$ is a trigonometric
polynomial of degree $N$, that the coefficients $x^{(k)}_{\n}$
are trigonometric polynomials with the stated properties
can be proved from \equ(2.6) by induction on $k$. \qed

\*

Then the functions $x^{(k)}(t)$ are well defined to all orders.
Before discussing the issue of convergence of the formal power series
defining such functions we look for a graphical representation
of the coefficients $x^{(k)}_{\n}$.

\*\*
\section(3,Graphical representation and tree formalism)

\0We start by giving some abstract definitions.

\*

\0{\bf Definition 1 (Trees).}
{\it A tree $\th$ is a graph, that is a connected set of points
and lines, with no cycle, such that all the lines are oriented toward
a unique point which has only one incident line.
Such a point is called the {\rm root} of the tree.
All the points in a tree except the root are denoted {\rm nodes}.
The line entering the root is called the {\rm root line}. 
The orientation of the lines in a tree induces a partial ordering 
relation between the nodes. We denote this relation by $\preceq$:
given two nodes $v$ and $w$, we shall write $w \preceq v$ every
time $v$ is along the path (of lines) which connects $w$ to the root.}

\*

Given a tree $\th$, we can identify the following subsets in $\th$.

\*

\0{\bf Definition 2 (Endpoints).}
{\it We call $E(\th)$ the set of {\rm endpoints} in $\th$,
that is the nodes which have no entering line. The endpoints can be
represented either as {\rm white bullets} or as {\rm black bullets}.
We call $E_W(\th)$ the set of white bullets and $E_{B}(\th)$
the set of black bullets. Of course $E_{W}(\th)\cup E_{B}(\th)=E(\th)$.
With each $v \in E_{W}(\th)$ we associate a {\rm mode} label $\n_v = 0$,
an {\rm order} label $k_{v} \in \ZZZ_{+}$ and a {\rm node factor}
$F_v = c_{k_{v}}$. With each $v \in E_{B}(\th)$ we associate
a {\rm mode} label $\n_{v}\in\ZZZ\setminus\{0\}$,
and a {\rm node factor} $F_{v} = f_{\n_{v}}$.}

\*

\0{\bf Definition 3 (Lines).}
{\it We denote with $L(\th)$ the set of {\rm lines} in $\th$.
Each line $\ell \in L(\th)$ leaves a point $v$ and enters another
one which we shall denote by $v'$. Since $\ell$ is uniquely
identified with $v$ (the point which $\ell$ leaves),
we may write $\ell = \ell_v$. 
With each line $\ell$ we associate a {\rm momentum} label $\n_{\ell}
\in \ZZZ$ and a {\rm propagator}
$$ g_{\ell} =
\cases{
1/(i\o\n_{\ell}) , & $\n_{\ell} \neq 0$, \cr
1 , & $\n_{\ell}=0$, \cr}
\Eq(3.1) $$
and we say that the momentum $\n_{\ell}$ flows through the line $\ell$.
The modes and the momenta are related as follows:
if $\ell = \ell_{v}$ one has
$$ \n_{\ell} = \sum_{i=1}^{s_{v}} \n_{\ell_{i}} =
\sum_{w \in E_{B}(\th) \atop w \preceq v} \n_{w} ,
\Eq(3.2) $$
where $\ell_{1},\ldots,\ell_{s_{v}}$ are the lines entering $v$.}

\*

\0{\bf Definition 4 (Vertices).}
{\it We denote by $V(\th)$ the set of {\rm vertices} in $\th$,
that is the set of points which have at least one entering line.
If $V(\th)\neq\emptyset$ we call the vertex $v_{0}$
connected to the root the {\rm last vertex} of the tree.
If $s_{v}$ denotes the number of lines entering $v$ call
$\max_{v\in V(\th)} s_{v}$ the {\rm branching number}.
One can have either $s_{v}=1$ or $s_{v}=2$.
We set $V_{s}(\th)=\{v\in V(\th) : s_{v}=s\}$ for $s=1,2$;
of course $V_{1}(\th)\cup V_{2}(\th) = V(\th)$. We define also
$V_{0}(\th)=\{v\in V(\th) : \n_{\ell_{v}}=0\}$; one has
$V_{0}(\th) \subset V_{2}(\th)$. We require that either
$V_{0}(\th)=\emptyset$ or $V_{0}(\th)=\{v_{0}\}$,
and that one can have $v\in V_{1}(\th)$ only if $\n_{\ell_{v}}\neq0$.
We associate with each vertex $v\in V(\th)$ a {\rm node factor}
$$ F_{v} = \cases{
-1 , & $s_{v} = 2$ and $v\notin V_{0}(\th)$ , \cr
-1/2c_{0} , & $s_{v} = 2$ and $v\in V_{0}(\th)$ , \cr
-(i\o\n_{\ell_{v}})^{2} , & $s_{v}=1$, \cr}
\Eq(3.3) $$
which is always well defined as $c_{0} \neq 0$.}

\*

We call {\it equivalent} two trees which can be transformed into
each other by continuously deforming the lines in such a way that
they do not cross each other.

Let $\TT_{k,\n}$ be the set of inequivalent trees $\th$ such that\\
(1) the number of vertices, the number of black bullets
and the order labels of the white bullets are such that we have
$$ \cases{
k_{1} + k_{2} + k_{3} = k , & if $\n\neq0$ , \cr
k_{1} + k_{2} + k_{3} = k + 1 , & if $\n=0$ , \cr}
\Eq(3.4) $$
if we set $k_{1} = |V(\th)|$, $k_{2}=|E_{B}(\th)|$
and $k_{3} = \sum_{v \in E_W(\th)} k_v$, and\\
(2) the momentum flowing through the root line is $\n$.

We refer to $\TT_{k,\n}$ as the {\it set of trees of order
$k$ and total momentum $\n$.}

With the above definitions the following result holds.

\*

\0{\bf Lemma 2.} {\it For all $k \ge 1$ and all $\n\neq0$ one has
$$ x^{(k)}_{\n} = \sum_{\th \in \TT_{k,\n}} \Val(\th) , \qquad
\Val(\th) = \left(\prod_{\ell \in L(\th)} g_{\ell} \right) 
\left(\prod_{v \in E(\th) \cup V(\th)} F_{v} \right) ,
\Eq(3.5) $$
where $\Val: \TT_{k,\n} \to \CCC$ is called the {\rm value} of the tree.
For $k\ge 2$ and $\n=0$ one has
$$ x^{(k)}_{0} \= c_{k} =
{\mathop{\sum}_{\th\in\TT_{k,0}}}^{*} \Val(\th) ,
\Eq(3.6) $$
where $*$ means that there are two lines entering the last
vertex $v_{0}$ of $\th_{0}$ and neither one exits from
an endpoint $v$ with order label $k_{v}=0$.}

\*

\0{\it Proof.}
We can represent graphically $x^{(k)}_{0}=c_{k}$ as in
Figure \graf(3.1)a, $x^{(1)}_{\n}$, $\n\ne0$, as in Figure \graf(3.1)b,
and, more generally, $x^{(k)}_{\n}$ as in Figure \graf(3.1)c.

\midinsert
\*
\insertplotttt{271pt}{19pt}{
\ins{-15pt}{15pt}{(a)}
\ins{35pt}{26pt}{$(k)$}
\ins{90pt}{15pt}{(b)}
\ins{143pt}{26pt}{$(1)$}
\ins{126pt}{6pt}{$\n$}
\ins{200pt}{15pt}{(c)}
\ins{255.5pt}{30.5pt}{$(k)$}
\ins{233pt}{6pt}{$\n$}
}
{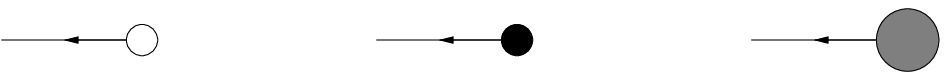}
\*
\line{\vtop{\line{\hskip1.5truecm\vbox{\advance\hsize by -3.1 truecm
\0{\nota
\figura(3.1).
Graphical representation of $x^{(k)}_{0}$, $x^{(1)}_{\n}$
and $x^{(k)}_{\n}$. For $\n=0$ the latter reduces
to the first graph, while for $k=1$ and $\n\ne0$ it reduces
to the second graph. In the first graph the momentum
is not shown as it is necessarily $\n=0$.}
} \hfill} }}
\*
\endinsert

Then the third equation in \equ(2.6) can be represented
graphically as in Figure \graf(3.2),
if we associate with the nodes and to the lines the node factors
and the propagators, respectively, according to 
the definitions \equ(3.1) and \equ(3.3).

\midinsert
\*
\insertplotttt{271pt}{59pt}{
\ins{40pt}{50pt}{$(k)$}
\ins{16pt}{26pt}{$\n$}
\ins{63pt}{31pt}{$=$}
\ins{89pt}{26pt}{$\n$}
\ins{141pt}{50pt}{$(k\!-\!1)$}
\ins{125pt}{26pt}{$\n$}
\ins{169pt}{36pt}{$+$}
\ins{197pt}{26pt}{$\n$}
\ins{255pt}{72pt}{$(k_{1})$}
\ins{255pt}{30pt}{$(k_{2})$}
\ins{237pt}{36pt}{$\n_{1}$}
\ins{237pt}{12pt}{$\n_{2}$}
}
{fishfig2}
\*
\line{\vtop{\line{\hskip1.5truecm\vbox{\advance\hsize by -3.1 truecm
\0{\nota
\figura(3.2).
Graphical representation of the third equation in \equ(2.6)
expressing the coefficient $x^{(k)}_{\n}$ for $k\ge2$ and $\n\neq0$
in terms of the coefficients  $x^{(k')}_{\n'}$ with $k'<k$.
In the last graph one has the constraints $k_{1}+k_{2}=k-1$
and $\n_{1}+\n_{2}=\n$.}
} \hfill} }}
\*
\endinsert

Analogously \equ(2.10) is represented graphically
as in Figure \graf(3.3), again if we use the graphical
representations in Figure \graf(3.1) and associate with the
lines and vertices the propagators \equ(3.1)
and the node factors \equ(3.3), respectively.

\midinsert
\*
\insertplotttt{163pt}{59pt}{
\ins{35pt}{47pt}{$(k)$}
\ins{52pt}{31pt}{$=$}
\ins{147pt}{71pt}{$(k_{1})$}
\ins{147pt}{31pt}{$(k_{2})$}
\ins{126pt}{35pt}{$\n_{1}$}
\ins{126pt}{15pt}{$\n_{2}$}
}
{fishfig3}
\*
\line{\vtop{\line{\hskip1.5truecm\vbox{\advance\hsize by -3.1 truecm
\0{\nota
\figura(3.3).
Graphical representation of the equation \equ(2.10)
expressing the coefficient $c_{k}$ for $k\ge2$ in terms of the
coefficients  $x^{(k')}_{\n'}$ with $k'<k$. Both $k_{1}$ and $k_{2}$
are strictly positive and $k_{1}+k_{2}=k$; moreover $\n_{1}+\n_{2}=0$.}
} \hfill} }}
\*
\endinsert

Note that in this way we represent graphically each
coefficient $x^{(k)}_{\n}$ in terms of other coefficients
$x^{(k')}_{\n'}$, with $k'<k$, so that we can apply iteratively
the graphical representation in Figure \graf(3.2) until only trees
whose endpoints represent either $x^{(1)}_{\n}$ with $\n\neq 0$
(black bullets) or $c_{k}$ are left (white bullets). This corresponds
exactly to the expressions in \equ(3.5) and \equ(3.6). \qed

\*

To get familiar with the graphic representation \equ(3.5) and \equ(3.6)
one should try to draw the trees which correspond to the first orders,
and check that the sum of the values obtained with the graphical
rules listed above gives exactly the same analytical expression
which can be deduced directly from \equ(2.6) and \equ(2.10).
 
For instance for $k=2$ we obtain for $x^{(2)}_{\n}$, $\n\ne0$,
the graphical representation in Figure \graf(3.4) and for
$c_{2}=x^{(2)}_{0}$ the graphical representation in Figure \graf(3.5).

\midinsert
\*
\insertplotttt{262pt}{49pt}{
\ins{39pt}{47pt}{$(2)$}
\ins{16pt}{21pt}{$\n$}
\ins{59pt}{25pt}{$=$}
\ins{85pt}{21pt}{$\n$}
\ins{143pt}{41pt}{$(1)$}
\ins{124pt}{21pt}{$\n$}
\ins{164pt}{27pt}{$+$}
\ins{193pt}{21pt}{$\n$}
\ins{251pt}{61pt}{$(0)$}
\ins{251pt}{22pt}{$(1)$}
\ins{234pt}{9pt}{$\n$}
}
{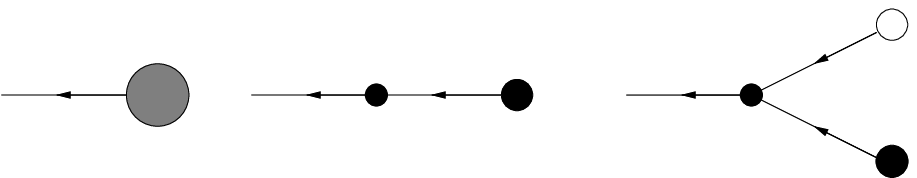}
\*
\line{\vtop{\line{\hskip1.5truecm\vbox{\advance\hsize by -3.1 truecm
\0{\nota
\figura(3.4).
Graphical representation of $x^{(2)}_{\n}$ for $\n\neq 0$.
The second contribution has to be counted twice,
because there is also a tree with the white and black bullets
exchanged, -- of course the latter has the same value.}
} \hfill} }}
\*
\endinsert

\midinsert
\*
\insertplotttt{154pt}{50pt}{
\ins{34pt}{43pt}{$(2)$}
\ins{52pt}{26pt}{$=$}
\ins{142pt}{63pt}{$(1)$}
\ins{142pt}{23pt}{$(1)$}
\ins{126pt}{30pt}{$\n_{1}$}
\ins{120pt}{8pt}{$-\n_{1}$}
}
{fishfig5}
\*
\line{\vtop{\line{\hskip1.5truecm\vbox{\advance\hsize by -3.1 truecm
\0{\nota
\figura(3.5).
Graphical representation of $c_{2}=x^{(2)}_{0}$. There is no
contribution with any white bullet carrying order label $k=0$
and $k=1$ because of the restriction in the sum appearing
in {\equ(3.6)} and of the fact that $c_{1}=0$, respectively.}
} \hfill} }}
\*
\endinsert

For $k=3$ we obtain for $x^{(3)}_{\n}$, $\n\ne0$, the graphical
representation in Figure \graf(3.6) and for $c_{3}=x^{(3)}_{0}$
the graphical representation in Figure \graf(3.7),
where we have explicitly used that $c_{1} = 0$.

\midinsert
\*
\insertplotttt{334pt}{194pt}{
\ins{39pt}{192pt}{$(3)$}
\ins{16pt}{166pt}{$\n$}
\ins{61pt}{172pt}{$=$}
\ins{90pt}{166pt}{$\n$}
\ins{125pt}{166pt}{$\n$}
\ins{180pt}{187pt}{$(1)$}
\ins{160pt}{166pt}{$\n$}
\ins{201pt}{174pt}{$+$}
\ins{233pt}{166pt}{$\n$}
\ins{267pt}{166pt}{$\n$}
\ins{324pt}{207pt}{$(0)$}
\ins{324pt}{167pt}{$(1)$}
\ins{307pt}{154pt}{$\n$}
\ins{61pt}{101pt}{$+$}
\ins{90pt}{94pt}{$\n$}
\ins{143pt}{135pt}{$(0)$}
\ins{125pt}{94pt}{$\n$}
\ins{180pt}{135pt}{$(0)$}
\ins{180pt}{95pt}{$(1)$}
\ins{164pt}{82pt}{$\n$}
\ins{202pt}{102pt}{$+$}
\ins{233pt}{94pt}{$\n$}
\ins{288pt}{135pt}{$(0)$}
\ins{267pt}{84pt}{$\n$}
\ins{324pt}{75pt}{$(1)$}
\ins{303pt}{67pt}{$\n$}
\ins{60pt}{28pt}{$+$}
\ins{90pt}{22pt}{$\n$}
\ins{145pt}{62pt}{$(1)$}
\ins{125pt}{32pt}{$\n_{1}$}
\ins{145pt}{22pt}{$(1)$}
\ins{125pt}{9pt}{$\n_{2}$}
}
{fishfig6}
\*
\line{\vtop{\line{\hskip1.5truecm\vbox{\advance\hsize by -3.1 truecm
\0{\nota
\figura(3.6).
Graphical representation of $x^{(3)}_{\n}$ for $\n\neq 0$.
The second and fourth contributions have to be counted twice,
while the third one has to be counted four times. There is no
contribution with any white bullet carrying the order label $k=1$
as $c_{1}=0$.}
} \hfill} }}
\*
\endinsert

\midinsert
\*
\insertplotttt{334pt}{70pt}{
\ins{35pt}{61pt}{$(3)$}
\ins{61pt}{44.5pt}{$=$}
\ins{90pt}{41pt}{$0$}
\ins{143pt}{82pt}{$(1)$}
\ins{118pt}{69pt}{$-\n$}
\ins{125pt}{39pt}{$\n$}
\ins{180pt}{82pt}{$(0)$}
\ins{180pt}{40pt}{$(1)$}
\ins{160pt}{29pt}{$\n$}
\ins{202pt}{39pt}{$+$}
\ins{233pt}{41pt}{$0$}
\ins{288pt}{82pt}{$(1)$}
\ins{267pt}{29pt}{$\n$}
\ins{260pt}{66pt}{$-\n$}
\ins{324pt}{22pt}{$(1)$}
\ins{303pt}{12pt}{$\n$}
}
{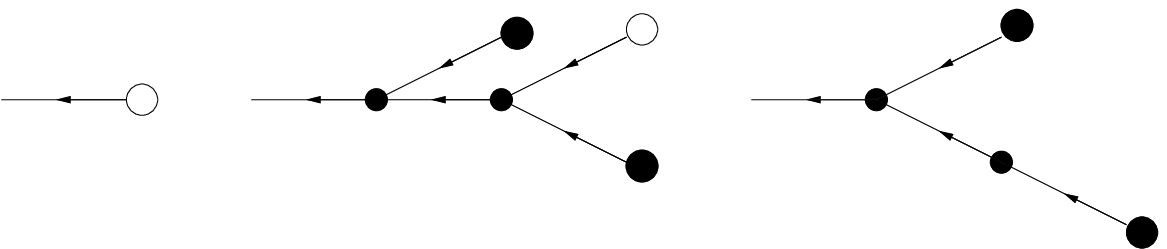}
\*
\line{\vtop{\line{\hskip1.5truecm\vbox{\advance\hsize by -3.1 truecm
\0{\nota
\figura(3.7).
Graphical representation of $c_{3}=x^{(3)}_{0}$.
The second contribution has to be counted twice,
while the first one has to be counted four times. There is no
contribution with any white bullet carrying the order label $k=1$
as $c_{1}=0$.}
} \hfill} }}
\*
\endinsert

And so on one can continue to higher orders. In general a
tree $\th\in\TT_{k,\n}$ looks like in Figure \graf(3.8),
where for simplicity no labels have been drawn other
than the order labels of the white bullets.
Note that each node can have only one or two entering lines,
while the endpoints have no entering line at all.
Moreover the momentum flowing through the line exiting a vertex $v$
is equal to the sum of the momenta flowing through the lines
entering $v$, according to \equ(3.2): this is a sort of conservation law.
The order of the tree is given by the number of vertices and
black bullets plus the sum of the order labels of the white bullets
minus the number of vertices in $V_{0}(\th)$. The latter
is just $|V_{0}(\th)|=0$ if $\th\in\TT_{k,\n}$, $\n\neq0$,
and $|V_{0}(\th)|=1$ if $\th\in\TT_{k,0}$.

\midinsert
\*
\insertplotttt{262pt}{86pt}{
\ins{68.6pt}{81pt}{$(k_{1})$}
\ins{176.6pt}{99pt}{$(k_{2})$}
\ins{248pt}{23pt}{$(k_{3})$}
}
{fishfig8}
\*
\line{\vtop{\line{\hskip1.5truecm\vbox{\advance\hsize by -3.1 truecm
\0{\nota
\figura(3.8).
Example of tree appearing in the graphical expansions {\equ(3.5)}
and {\equ(3.6)}. The number of lines entering any vertex $v$ can be
only either $s_{v}=1$ or $s_{v}=2$, while no line enters the endpoints.
The order of the tree is given by the number of elements
in $V(\th)\setminus V_{0}(\th)$ plus the number of elements
in $E_{B}(\th)$ plus the sum of the order labels of the white bullets.
Then, if $k_{1}$, $k_{2}$ and $k_{3}$ are the order labels of the
white bullets in the figure the order of the tree
is $k=k_{1}+k_{2}+k_{3}+10$ if $v_{0}\notin V_{0}(\th)$ and
$k=k_{1}+k_{2}+k_{3}+9$ if $v_{0}\in V_{0}(\th)$.
In the latter case one must have $k_{1}>0$
because of the constraint in the sum appearing in \equ(3.6).}
} \hfill} }}
\*
\endinsert

If a vertex $v$ has $s_{v}=1$, that is it has
only one entering line $\ell$, the latter can not come out from
a white bullet. Indeed if this occurs one should have
$\n_{\ell_{v}}=\n_{\ell}=0$, hence $F_{v}=0$ by \equ(3.3),
so that the value of the tree containing such a vertex is zero.

Given a tree as in Figure \graf(3.8) we can represent
each white bullet according to the graphic representation
in Figure \graf(3.3), corresponding to the analytic formula \equ(3.6),
and expand again the two contributions $x^{(k_{1})}_{\n_{1}}$
and $x^{(k_{2})}_{\n_{2}}$ as sums of trees,
and so on, iteratively, until the only white bullets which are
left are the ones with order label $k=0$. In this way we obtain a new
graphic representation where the trees still look like
in Figure \graf(3.8), but now there are a few differences:\\
(1) all the white bullets $v \in E_{W}(\th)$
have order labels $k_{v}=0$, and\\
(2) there can be lines $\ell\in L(\th)$ with momentum
$\n_{\ell}=0$ which come out from vertices, that is
$V_{0}(\th)$ can contain more than one or none element.

Note that only lines coming out either from nodes in
$V_{0}(\th) \subset V_{2}(\th)$ or
from white bullets have vanishing momentum.

The order of the tree is then given by
the number of elements of $V(\th)\cup E_{B}(\th)$ minus the number
of elements of $V_{0}(\th)$, that is
$k=|V(\th)|+|E_{B}(\th)|-|V_{0}(\th)|$.
Of course $v_{0}\in V_{0}(\th)$ if and only if the momentum of the
root line is vanishing, that is $\th\in\TT_{k,0}$ for some $k\ge2$.
It is important to stress that no line entering a vertex
$v \in V_{0}(\th)$ can come out from a white bullet
(which now has necessarily an order label $0$),
because this would be against the constraint in the sum \equ(3.6).
This means that if two lines carrying zero momentum enter the
same vertex $v$ (so that $v\in V_{0}(\th)$ according to \equ(3.2)),
then none of them can exit a white bullet.

But up to these minor differences a tree representation like
in \equ(3.5) and \equ(3.6) still holds. The advantage of these
modified rules is that now the tree values are expressed
no longer in terms of constants $c_{k}$ to be determined,
but only in terms of $c_{0}$ which is known. A tree drawn
according these new rules is represented as in Figure \graf(3.8)
with $k_{1}=k_{2}=k_{3}=0$ (and in particular a tree of this
kind can contribute only to $x^{(k)}_{\n}$ with $\n\neq0$).
Note that we could avoid drawing the order labels associated
with the endpoints, as they are uniquely determined as $k=0$
for the white bullets and $k=1$ for the black bullets.
Of course, with respect to the caption of that Figure,
now the order $k$ is given by the number of elements in $V(\th)$
plus the number of elements in $E_{B}(\th)$ minus
the number of elements in $V_{0}(\th)$.

\*\*
\section(4,Formal solutions)

\0The sum over the trees in \equ(3.5) and \equ(3.6), with the new
definition of the set $\TT_{k,\n}$ given at the end of Section \secc(3),
can be performed by summing over all possible ``tree shapes''
(that is trees without labels or {\it unlabelled trees}), and,
for a fixed shape, over all possible assignments of mode labels.
In the case of a trigonometric polynomial of degree $N$
the latter can be bounded by $(2N)^{|E(\th)|}$, because
each endpoint $v$ can have either a mode label $\n_{v}\neq 0$, with
$|\n_{v}|\le N$, or the mode label $\n_{v}=0$,
while the case of analytic functions
(or even to obtain bounds which are uniform in $N$)
has to be discussed a little more carefully.
The number of unlabelled trees with $P$ nodes
(vertices and endpoints) can be bounded by $2^{2P}$.

Recall that $V_{s}(\th)$ denotes the set of vertices $v$
such that $s_{v}=s$; of course $V_{1}(\th)\cup V_{2}(\th)=V(\th)$,
and $V_{0}(\th)\subset V_{2}(\th)$.
Analogously we can set
$$ \eqalign{
L_{0}(\th) & = \left\{ \ell \in L(\th) : n_{\ell}=0 , \right\} , \cr
L_{1}(\th) & = \left\{ \ell \in L(\th) : \ell=\ell_{v} ,
v \in V_{1}(\th) \right\} , \cr
L_{2}(\th) & = L(\th) \setminus \big( L_{0}(\th) \cup L_{1}(\th)
\big) , \cr}
\Eq(4.1) $$
with the splitting made in such a way that one has
$$ \eqalign{
& \Big| \prod_{v \in V_{1}(\th) } F_{v} \Big| \,
\Big| \prod_{\ell \in L_{1}(\th) } g_{\ell} \Big|
\le \prod_{\ell \in L_{1}(\th) } \left| \o \n_{\ell} \right| ,
\qquad \Big| \prod_{\ell \in L_{2}(\th) } g_{\ell} \Big|
\le \prod_{\ell \in L_{2}(\th)}
{1 \over \left| \o \n_{\ell} \right| } , \cr
& \Big| \prod_{v \in V_{0}(\th) } F_{v} \Big|
\le \left( {1 \over 2c_{0}} \right)^{|V_{0}(\th)|} , \qquad
\Big| \prod_{v \in E_{W}(\th) } F_{v} \Big|
\le c_{0}^{|E_{W}(\th)|} , \cr
&
\Big| \prod_{v \in E_{B}(\th) } F_{v} \Big|
\le F^{|E_{B}(\th)|}
\prod_{v \in E_{B}(\th) } {\rm e}^{-\x|\n_{v}|} , \cr}
\Eq(4.2) $$
where for each line $\ell$ one
has $|\n_{\ell}| \le \sum_{v\in E_{W}(\th)} |\n_{v}|$.

The following result is useful when
looking for bounds on the tree values.

\*

\0{\bf Lemma 3.} {\it Given a tree $\th$
with branching number $s$ one has $|E(\th)| \le (s-1)|V(\th)|+1$.
If $k$ denotes the order of the tree $\th$, that is
$|V(\th)|-|V_{0}(\th)|+|E_{B}(\th)|=k$, one has
the identity $|L_{1}(\th)|+|L_{2}(\th)|=k$, and
the bounds $|V_{1}(\th)|\le k$, $|V_{0}(\th)| \le k-1$,
$|E(\th)|\le k$ and $|E(\th)|+|V(\th)| \le 2k - 1$.}

\*

\0{\it Proof.} It is a standard result on trees that one has
$\sum_{v\in V(\th)}(s_{v}-1)=|E(\th)|-1$, so that the first
bound follows. The bounds on $|V_{1}(\th)|$, $|V_{0}(\th)|$, $|E(\th)|$
and $|E(\th)|+|V(\th)|$ can be easily proved by induction, while
the identity $|L_{1}(\th)|+|L_{2}(\th)|=k$ follows from
the observation that all lines in $L_{1}(\th)$ and $L_{2}(\th)$
come out either from vertices or from black bullets,
and they have non-vanishing momentum. \qed

\*

Hence the number of lines in $L_{1}(\th)$ is bounded by $k$,
so that in \equ(4.2) we can bound
$$ \eqalign{
\Big( \prod_{\ell \in L_{1}(\th)} |\o\n_{\ell}| \Big)
\Big( \prod_{v\in E_{B}(\th)} F \, {\rm e}^{-\x|\n_{v}|} \Big)
& \le \Big( \prod_{v\in E_{B}(\th)} F \, {\rm e}^{-\x|\n_{v}|/2} \Big)
\Big( \prod_{\ell \in L_{1}(\th)} {\rm e}^{-\x|\n_{\ell}|/2k}
|\o\n_{\ell}| \Big) \cr
& \le \Big( \prod_{v\in E_{B}(\th)} F \, {\rm e}^{-\x|\n_{v}|/2} \Big)
\left( {2 k |\o| \over \x} \right)^{k} , \cr}
\Eq(4.3) $$
and in the second line the product can be used to perform the sum
over the Fourier labels -- and this gives a factor $F^{k}B_{2}^{k}$,
with $B_{2}=2{\rm e}^{-\x/2}(1-{\rm e}^{-\x/2})^{-1}$, --
while the last factor is bounded by $A_{1}B_{1}^{k}k!$,
for some constants $A_{1}$ and $B_{1}$.

We can bound the value of a tree $\th$ by using
the bounds \equ(4.2) and \equ(4.3), and Lemma 3. If we define
$$ \e_{1}^{-1} = \max\{ B_{1}, |\o|^{-1}\} \,
\max\{c_{0},F B_{2}\} \, \max\{1,(2c_{0})^{-1}\} ,
\Eq(4.4) $$
with $c_{0}=\sqrt{\a}$, and take into account that the number
of unlabelled trees in $\TT_{k,\n}$ is bounded by $2^{2k-1}$
(because each tree in $\TT_{k,\n}$ has at most $2k-1$ nodes), then

$$ \left| x^{(k)}_{\n} \right| \le A_{1} \e_{2}^{-k}k! , \qquad
\left| x^{(k)}(t) \right| \le A_{1} \e_{2}^{-k}k! ,
\Eq(4.5) $$
where we have set $\e_{2}=\e_{1}2^{-2}$.

A bound like \equ(4.5) is obtained also in the case of
forcing terms which are trigonometric polynomials, because in general
we can bound the factors $|\o\n_{\ell}|$ in \equ(4.2) only with $kN$
(see Lemma 1), and this produces an overall bound proportional to $k!$
Note that in that case the bound $B_{2}$, arising from
the sum over the Fourier labels, can be replaced with
a factor $2N$, and $B_{1}$ can be replaced with $|\o|N$.

Then we have proved the following result.

\*

\0{\bf Proposition 1.} {\it Given the equation \equ(2.1) with $f$
as in \equ(1.4), there is only one periodic solution
in the form of a formal power series, and the corresponding period is
the same period $2\p/\o$ of the forcing term. The coefficients of
such a formal power series satisfy the bounds \equ(4.5).}

\*

One could ask if the factorials arising in the bounds are only
a technical problem, or else they are a symptom
that the series really diverges. To order $k$ one can easily provide
examples of trees which grow like factorials;
see for instance the tree represented in Figure \graf(4.1),
where there are $k-1$ vertices with only one entering line.
Then the corresponding value is
$$ \Val(\th) = (i\o\n)^{2(k-1)} {1 \over (i\o\n)^{k}} f_{\n} =
(i\o\n)^{k-2} f_{\n} ,
\Eq(4.6) $$
which behaves as $k!$ for large $k$. Furthermore
it is unlikely that there are cancellations with the values
of other trees because the value of any other tree $\th\in\TT_{k,\n}$
can be proportional at most to $(i\o\n)^{p}$, with $p<k-2$
(strictly). Hence we expect that the coefficients
$u^{(k)}_{\n}$, even if well defined to all orders,
grow like factorials, so preventing the convergence of the series.

\midinsert
\*
\insertplotttt{298pt}{10pt}{
\ins{13pt}{2pt}{$\n$}
\ins{51pt}{2pt}{$\n$}
\ins{87pt}{2pt}{$\n$}
\ins{123pt}{2pt}{$\n$}
\ins{161pt}{2pt}{$\n$}
\ins{198pt}{2pt}{$\n$}
\ins{232pt}{2pt}{$\n$}
\ins{270pt}{2pt}{$\n$}
\ins{288pt}{21pt}{$(1)$}
}
{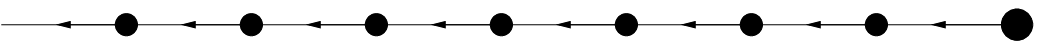}
\*
\line{\vtop{\line{\hskip1.5truecm\vbox{\advance\hsize by -3.1 truecm
\0{\nota
\figura(4.1).
Example of tree whose value grows as a factorial.
If $k$ is the order of the tree (hence there are $k-1$ vertices and $1$
black bullet), then the value of the tree is given in \equ(4.7).}
} \hfill} }}
\*
\endinsert

The lack of analyticity is further supported by the following fact.
If we consider \equ(2.1) without the quadratic term and with $\a=0$,
that is
$$ \e \ddot x + \dot x = f(\o t) , \qquad f_{0}=0 ,
\Eq(4.7) $$
in Fourier space, we find $x_{0}=0$ and $i\o\n (1+i\e\o\n) x_{\n} =
f_{\n}$ for $\n\neq0$. Hence the equation is trivially solvable,
and it gives
$$ x(t) = \sum_{\n\ne0} {f_{\n} \over i\o\n (1+i\e\o\n)}
{\rm e}^{i\o\n t} .
\Eq(4.8) $$
Of course the solution $x(t)$ of the linear equation
is not analytic in $\e$ (in a neighbourhood of the origin)
when $f$ is an analytic function containing all the harmonics,
as each point $\e=i/\o\n$ represents a singularity point for $x(t)$,
and such points accumulate to the origin as $\n\to\io$.
Then it is likely that also when the quadratic terms
are taken into account the solution can not be analytic.
Therefore giving a meaning to the perturbation series
requires some more work, that we are going to discuss next.

An important remark is that for any $k\ge 1$
there is no tree whose value can be bounded worse than proportionally
to a factorial, as the estimates \equ(4.5) show:
indeed they have been obtained by bounding separately
the value of each single tree. This observation
will play an important role in the forthcoming analysis.

\*\*
\section(5,Periodic forcing terms)

\0To deal completely with the case of analytic functions and prove
existence of the periodic solution, we have to modify the
graphical expansion envisaged in the previous sections.

Let us come back to the equation \equ(2.1), and write it in
Fourier space. For $\n\ne0$ and denoting with
$x_{\n}$ the $\n$-th Fourier coefficient, we obtain
$$ \e (i\o\n)^{2} x_{\n} + i\o\n \, x_{\n} + 
\e \sum_{\n_{1}+\n_{2}=\n} x_{\n_{1}} x_{\n_{2}} = \e f_{\n} ,
\Eq(5.1) $$
provided that for $\n=0$ we have
$$ \sum_{\n_{1}+\n_{2}=0} x_{\n_{1}} x_{\n_{2}} = 0 .
\Eq(5.2) $$
Let us rewrite \equ(5.1) as
$$ \e (i\o\n)^{2} x_{\n} + i\o\n \, x_{\n} + \m \e
\sum_{\n_{1}+\n_{2}=\n} x_{\n_{1}} x_{\n_{2}} = \m \e f_{\n} ,
\Eq(5.3) $$
and look for a solution $x(t)$ which is analytic in $\m$,
which suggests us to write
$$ x(t) = \sum_{k=0}^{\io} \m^{k} x^{[k]}(t) .
\Eq(5.4) $$
Of course we want that at the end the value $\m=1$ is inside
the analyticity domain. Note also that now $x^{[k]}$,
the coefficient to order $k$, has a different meaning with respect
to the previous expansion \equ(1.4) in powers of $\e$,
and for this reason with use a different symbol to denote it.
We shall call {\it resummed series} the series \equ(5.4),
because the coefficients $x^{[k]}(t)$ depend on $\e$,
and are given by the sum of infinitely many terms
of the formal series \equ(2.5).

Again for $k=0$ we have to take $x^{[0]}_{\n}=0$ for $\n\ne0$
and fix $c_{0}\=x^{[0]}_{0}=\sqrt{\a}$, with $\a\=f_{0}$.

To order $k\ge 1$ (in $\m$) we obtain for $\n\ne0$
$$ i \o \n \left( 1 + i \e \o \n \right) x_{\n}^{[k]} =
\e f_{\n} \d_{k,1} - \e \sum_{k_{1}+k_{2}=k-1}
\sum_{\n_{1}+\n_{2}=\n} x^{[k_{1}]}_{\n_{1}} x^{[k_{2}]}_{\n_{2}} ,
\Eq(5.5) $$
while for $\n=0$ we require
$$ \sum_{k_{1}+k_{2}=k} \sum_{\n_{1}+\n_{2}=\n}
x^{[k_{1}]}_{\n_{1}} x^{[k_{2}]}_{\n_{2}} = 0 .
\Eq(5.6) $$
By setting $c_{k}=x^{[k]}_{0}$ the latter equation can
be written as (cf. \equ(2.10))
$$ c_{1} = 0 , \qquad\qquad
c_{k} = -{1 \over 2 c_{0}} \sum_{k'=1}^{k-1}
\sum_{\n\in\zzz} x^{[k-k']}_{\n} x^{[k']}_{-\n} , \qquad k \ge 2 .
\Eq(5.7) $$

Then we can proceed as in Section \secc(2), with some slight changes
that we now explain. First of all note that \equ(5.5) gives
for $\n\neq0$
$$ \eqalign{
x^{[0]}_{\n} & = 0 , \cr
x^{[1]}_{\n} & = {\e f_{\n} \over i\o\n ( 1+i\e\o\n)} , \cr
x^{[k]}_{\n} & = - {\e \over i\o\n ( 1+i\e\o\n)}
\sum_{k_{1}+k_{2}=k-1}
\sum_{\n_{1}+\n_{2}=\n} x^{[k_{1}]}_{\n_{1}} x^{[k_{2}]}_{\n_{2}} ,
\qquad k\ge2 . \cr}
\Eq(5.8) $$
Then the graphical representations of $x^{[k]}_{0}$, $x^{[1]}_{\n}$
and $x^{[k]}_{\n}$ are as in the previous case, with the only change
in the representation of the order labels (because of the square
brackets instead of the parentheses); see Figure \graf(5.1).

\midinsert
\*
\insertplotttt{271pt}{19pt}{
\ins{-15pt}{15pt}{(a)}
\ins{35pt}{26pt}{$[k]$}
\ins{90pt}{15pt}{(b)}
\ins{145pt}{26pt}{$[1]$}
\ins{126pt}{6pt}{$\n$}
\ins{200pt}{15pt}{(c)}
\ins{257.5pt}{30.5pt}{$[k]$}
\ins{233pt}{6pt}{$\n$}
}
{fishfig1}
\*
\line{\vtop{\line{\hskip1.5truecm\vbox{\advance\hsize by -3.1 truecm
\0{\nota
\figura(5.1).
Graphical representation of $x^{[k]}_{0}$, $x^{[1]}_{\n}$
and $x^{[k]}_{\n}$. For $\n=0$ the latter reduces
to the first graph, while for $k=1$ and $\n\ne0$ it reduces
to the second graph. In the first graph the momentum
is not showed as it is necessarily $\n=0$.}
} \hfill} }}
\*
\endinsert

On the contrary the graphical representation of the third equation
in \equ(5.8) is as in Figure \graf(5.2).

\midinsert
\*
\insertplotttt{169pt}{59pt}{
\ins{41.pt}{51pt}{$[k]$}
\ins{16.7pt}{26.7pt}{$\n$}
\ins{61.7pt}{30.8pt}{$=$}
\ins{89.2pt}{26.7pt}{$\n$}
\ins{150pt}{70.5pt}{$[k_{1}]$}
\ins{150pt}{30.5pt}{$[k_{2}]$}
\ins{125.8pt}{35pt}{$\n_{1}$}
\ins{125.8pt}{12.5pt}{$\n_{2}$}
}
{fishfig10}
\*
\line{\vtop{\line{\hskip1.5truecm\vbox{\advance\hsize by -3.1 truecm
\0{\nota
\figura(5.2).
Graphical representation of the second equation in \equ(5.8)
expressing the coefficient $x^{[k]}_{\n}$ for $k\ge 2$ and $\n\neq0$
in terms of the coefficients  $x^{[k']}_{\n'}$ with $k'<k$.
In the last graph one has the constraints $k_{1}+k_{2}=k-1$
and $\n_{1}+\n_{2}=\n$.}
} \hfill} }}
\*
\endinsert

At the end we obtain a tree expansion where the trees differ from the
previous ones as they contain no vertex with only one entering line. 
With the previous notations this means that $L_{1}(\th)=\emptyset$
and $V_{1}(\th)=\emptyset$, hence $V(\th)=V_{2}(\th)$.
Moreover also the propagators and the node factors of the vertices
are different, as \equ(3.1) and \equ(3.2) have to be replaced with
$$ g_{\ell} =
\cases{
1/((i\o\n_{\ell})(1+i\e\o\n_{\ell})) , & $\n_{\ell} \neq 0$, \cr
1 , & $\n_{\ell}=0$, \cr}
\Eq(5.9) $$
and, respectively,
$$ F_{v} = \cases{
-\e , & $v \notin V_{0}(\th)$ , \cr
-1/2c_{0} , & $v \in V_{0}(\th)$ , \cr}
\Eq(5.10) $$
and we recall once more that only vertices $v$ with $s_{v}=2$
are allowed. Finally, the node factors associated to the
endpoints are $F_{v}=c_{k_{v}}$ if $v$ is a white bullet
and $F_{v}=\e f_{\n_{v}}$ if $v$ is a black bullet.

As in Section \secc(3) we can envisage an expansion in which
all white bullets $v$ have $k_{v}=0$ (simply by expanding
iteratively in trees the white bullets of higher order).
A tree appearing in this new expansion is
represented in Figure \graf(5.3).

\midinsert
\*
\insertplotttt{262pt}{86pt}{
}
{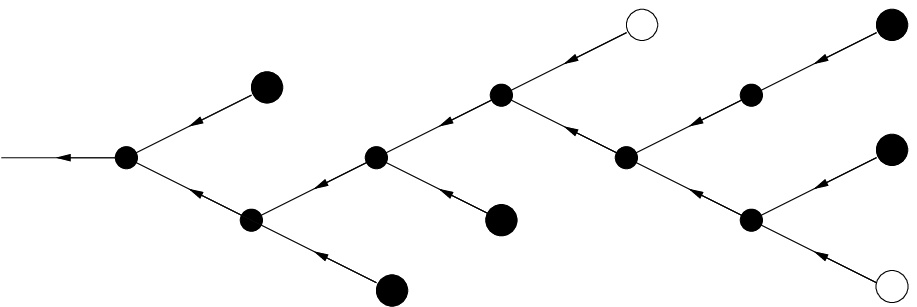}
\*
\line{\vtop{\line{\hskip1.5truecm\vbox{\advance\hsize by -3.1 truecm
\0{\nota
\figura(5.3).
Example of tree appearing in the new graphical expansion.
The number of lines entering any vertex $v$ can be only
$s_{v}=2$. The order of the tree is given by $|B(\th)|-|V_{0}(\th)|$.
All the white bullets have order labels $0$,
as well as all the black bullets carry a label $(1)$;
hence we can avoid drawing explicitly such labels.}
} \hfill} }}
\*
\endinsert

With the notations \equ(4.1), we obtain the bounds
$$ \eqalign{
& \Big| \prod_{v \in V(\th) \setminus V_{0}(\th)} F_{v} \Big|
\le |\e|^{|V(\th)|} ,
\qquad \Big| \prod_{\ell \in L(\th) } g_{\ell} \Big|
\le \prod_{\ell \in L_{2}(\th) }
{1 \over \left| \o \n_{\ell} \right| \,
\left| 1 + i \e \o \n_{\ell} \right|} , \cr
& \Big| \prod_{v \in V_{0}(\th) } F_{v} \Big|
\le \left( {1 \over 2c_{0}} \right)^{|V_{0}(\th)|} , \qquad
\Big| \prod_{v \in E_{W}(\th) } F_{v} \Big|
\le c_{0}^{|E_{W}(\th)|} , \cr
& \Big| \prod_{v \in E_{B}(\th) } F_{v} \Big|
\le F^{|E_{B}(\th)|}
\prod_{v \in E_{B}(\th) } {\rm e}^{-\x|\n_{v}|} , \cr}
\Eq(5.11) $$
where we have used again the bound $|f_{\n}|\le F\,{\rm e}^{-\x|\n|}$, for
suitable (strictly) positive constants $F$ and $\x$, which follows
from the analyticity assumption on $f$.

For real $\e$ we can bound each propagator by
$$ \left| g_{\ell} \right| \le {1 \over |\o \n_{\ell}|}
\le {1 \over |\o|} ,
\Eq(5.12) $$
so that the value of any tree $\th\in\TT_{k,\n}$ can be bounded by
$$ \left| \Val(\th) \right| \le |\e|^{k}
|\o|^{-k} (\max\{c_{0},F\})^{k} (\max\{1,1/2c_{0}\})^{k}
\prod_{v \in E_{B}(\th)} {\rm e}^{-\x|\n_{v}|} ,
\Eq(5.13) $$
where we have used again Lemma 3. If we write
$$ \prod_{v \in E_{B}(\th)} {\rm e}^{-\x|\n_{v}|} \le {\rm e}^{-\x|\n|/2}
\left( \prod_{v \in E_{B}(\th)} {\rm e}^{-\x|\n_{v}|/2} \right) ,
\Eq(5.14) $$
we can proceed as in ection \secc(4): we use the last product to perform
the sum over the Fourier labels, which gives a factor $B_{2}^{k}$,
whereas the sum over the unlabelled trees gives a factor $2^{2k-1}$.
At the end we obtain
$$ \left| x^{(k)}_{\n} \right| \le \m_{2}^{-k} , \qquad
\left| x^{(k)}(t) \right| \le \m_{2}^{-k} ,
\Eq(5.15) $$
where we have set $\m_{2}^{-1}=4\,|\o|^{-1} \max\{1,1/2c_{0}\}
\max\{F B_{2} ,c_{0}\} |\e|$. Hence the radius of convergence $\m_{0}$
of the series expansion \equ(5.4) is bounded as
$\m_{0} \ge \m_{2} = O(1/|\e|)$, so that for $\e$ small enough,
say $|\e|<\e_{3}= (4\,|\o|^{-1} \max\{1,1/2c_{0}\}
\max\{F B_{2} ,c_{0}\})^{-1}$, the value $\m=1$
is inside the analyticity domain.

We can summarize the results found so far as follows.

\*

\0{\bf Theorem 1.} {\it Given the equation (2.1) with
$f$ analytic, there exists $\e_{0}>0$ such that for all real $\e$
with $|\e|<\e_{0}$ there is only one periodic solution
which admits a formal expansion in powers of $\e$,
and the corresponding period is the same period $2\p/\o$
of the forcing term. An explicit bound is $\e_{0} \ge \e_{3}=O(\o)$.}

\*

Note that if $\o$ is very large then very large values
of $\e$ are allowed.

We can investigate further the regularity properties in $\e$
of the periodic solution found in Theorem 1,
and see what happens for complex values of $\e$.

We need the following preliminary result (see Figure \graf(5.4)a
for the region $\CC_{R}$).

\midinsert
\*
\insertplotttt{325pt}{109pt}{
\ins{-15pt}{60pt}{(a)}
\ins{165pt}{60pt}{(b)}
\ins{135pt}{50pt}{${\rm Re}\,\e$}
\ins{50pt}{110pt}{${\rm Im}\,\e$}
\ins{317pt}{50pt}{${\rm Re}\,\e$}
\ins{195pt}{110pt}{${\rm Im}\,\e$}
}
{fishfig12}
\*
\line{\vtop{\line{\hskip1.5truecm\vbox{\advance\hsize by -3.1 truecm
\0{\nota
\figura(5.4).
Region $\CC_{R}$ in the complex $\e$-plane (a) and strip-like region
of analyticity $\SS_{B}$ of the Borel transform (b).
The region $\CC_{R}$ is the union of two discs of radius $R/2$ and
centers $(\pm R/2,0)$.}
} \hfill} }}
\*
\endinsert

\*

\0{\bf Lemma 4.} {\it Given $\o>0$ and $0<R<1/4\o$ let $\CC_{R}$
be the pair of discs $\CC_{R}=\{\e : |{\rm Re}\,\e^{-1}|>R^{-1}\}$.
For all $\e\in \CC_{R}$ and all $\n\in\ZZZ \setminus \{0\}$
one has $|i\o\n(1+i\e\o\n)| \ge \o/2$.}

\*

\0{\it Proof.} Write $\e=a+ib$ and $x=\o\n$, so that one has $|i\o\n
( 1+i\e\o\n)|=|x|\sqrt{(1-bx)^{2}+(ax)^{2}}\= F(x)$.
If $\e\in \CC_{R}$ one has $|a| \ge b^{2}/2R$. Fix $0<A<1$.
If $|1-bx| \le A$ then $\sqrt{(1-bx)^{2}+(ax)^{2}}
\ge |ax| \ge b^{2}|x|/2R \ge |b|(1-A)/2R$, so that
$F(x) \ge (1-A)^{2}/2R$.
If $|1-bx| \ge A$ then $\sqrt{(1-bx)^{2}+(ax)^{2}}
\ge A$, hence $F(x) \ge A|x| \ge \o A$. Then choose
$A=1-\sqrt{\o R}\ge 1/2$; this gives $F(x) \ge \o/2$. \qed

\*

Now fix $0<R < \lis R \= \e_{3}$ so small
that $|\o| R<1/4$, and consider the corresponding domain $\CC_{R}$.
We can apply Lemma 4 and deduce that any propagator
$g_{\ell}$ is bounded by $|g_{\ell}| \le 2/|\o|$ for all $\e\in \CC_{R}$.

This allows us to obtain the following result.

\*

\0{\bf Proposition 2.} {\it There exists $R>0$ small enough such that
in the domain $\CC_{R}$ one has the asymptotic expansion
$$ x(t) = \sum_{k=0}^{N-1} \e^{k} x^{(k)}(t) + \RRRR_{N}(\e) ,
\qquad |\RRRR_{N}(\e)| \le A B^{N} N! |\e|^{N} ,
\Eq(5.16) $$
where the constants $A$ and $B$ are uniform in $N$ and in $\e$.}

\*

\0{\it Proof.} Write $x(t)$ as $x(t)=x_{N}(t)+\RRRR_{N}(t)$, where
$x_{N}(t)$ is given by the sum of the first $N-1$ orders of the formal
power series expansion of the solution $x(t)$ as in \equ(5.16).
For $\e\in \CC_{R}$ the function \equ(5.4) with $\m=1$
is $C^{\io}$ in $\e$, hence we can estimate $\RRRR_{N}(\e)$
with a bound on the $N$-th derivative of $x(t)$ in $\CC_{R}$,
and this gives the bound in \equ(5.16). \qed

\*

Of course the constants $A$ and $B$ in \equ(5.16) are explicitly
computable; in particular one finds $B=O(\e_{3}^{-1})$.

Then we are under the assumptions where Nevanlinna's theorem
\cita{N} (see also Ref.~\cita{S}) can be applied,
and we obtain that the function
$$ \BBBB(t;\e) = \sum_{k=0}^{\io} {1 \over k!} \e^{k} x^{(k)}(t)
\Eq(5.17) $$
converges for $|\e|<B$ (with $B$ given in \equ(5.16)) and has an
analytic continuation to $\SS_{B}= \{ \e: {\rm dist}(\e,\RRR_{+})<B\}$
(see Figure \graf(5.4)b), satisfying for some constant $K$
the bound $| \BBBB(t;\e)| \le K {\rm e}^{|\e|/R}$
uniformly in every $\SS_{B'}$ with $B'<B$. The function
$x(t)$ can be represented as the absolutely convergent integral
$$ x(t) = {1\over \e} \int_{0}^{\io} e^{-s/\e} \BBBB(t;s) \, {\rm d}s
\Eq(5.18) $$
for all $\e\in \CC_{R}$, and this property can be stated by saying
that $x(t)$ is Borel-summable (in $\e$) and $\BBBB(t;\e)$ is its
Borel transform \cita{H}. This implies that the function given by
the summation procedure described in Theorem 1 is unique.
Therefore we have obtained the following result,
which strengthens Theorem 1.

\*

\0{\bf Theorem 2}. {\it The solution given by Theorem 1 is
Borel-summable at the origin.}

\*

Note that Watson's theorem can not be invoked to obtain this result
because the singularities are along the imaginary axis.

In particular if $f(\o t)=\a+\b \sin t$ then there is a periodic solution
$x(t) = a + \e \b \cos t + O(\e^{2})$, with $a = \sqrt{\a}$, which has
period $2\p$ and moves around the fixed point $(x,\dot x)=(a,0)$, and
close to it within $O(\e)$. No other periodic solution analytic in $\e$
can exist.

We conclude this section with two remarks. The summation
criterion envisaged in this Section is reminiscent of that
used (in a more difficult situation) in Ref.~\cita{GG1}
for hyperbolic lower-dimensional tori.
However in that case we have not been able to prove
Borel summability because to order $k$ the bounds are
like $(k!)^{\a}$ for some $\a>1$. Neither extensions to
Watson's theorem \cita{H} analogous to Nevanlinna-Sokal's result
(as those developed in Ref.~\cita{CGM})
can be used because the exponent $\a$ is too large.
We shall find a very similar situation in next section.

The lack of analyticity in $\e$ in a neighbourhood
of the origin is due to the accumulation
of singularity points along the imaginary axis in the
complex $\e$-plane (where the quantity $1+i\e\o\n$ vanishes
for $\n\in\ZZZ$). The analyticity domain is tangential
to the imaginary axis, and this allows us to apply
Nevanlinna's theorem. We find that this situation
has some analogies with a different problem,
the analyticity properties of rescaled versions
of some dynamical systems, such as Siegel's problem \cita{BMS}
-- and its linearization as considered in Ref.~\cita{MS} --,
the standard map \cita{BG1} and generalized standard maps \cita{BG2},
for complex rotation numbers tending to rational values
in the complex plane. In those cases, however,
only non-tangential limits could be considered.
Of course the situation is slightly more complicated there,
because the set of accumulating singularity points
is dense -- and not only numerable as in the present case.

\*\*
\section(6,Quasi-periodic forcing terms)

\0In the case of analytic quasi-periodic forcing terms, we shall assume
a Diophantine condition on the rotation vector $\oo$, that is
$$ \left| \oo\cdot\nn \right| \ge C_{0} |\nn|^{-\t} \qquad
\forall \nn\in\ZZZ^{d} \setminus\{\V0\} ,
\Eq(6.1) $$
where $|\nn|=|\nn|_{1}\=|\n_{1}|+\ldots+|\n_{d}|$,
and $C_{0}$ and $\t$ are positive constants.
We need $\t\ge d-1$ in order to have
a non-void set of vectors satisfying the condition \equ(6.1),
and $\t>d-1$ in order to have a full measure set of such vectors.
For simplicity (and without loss of generality) we can assume
$C_{0}<\g/2$, with $\g=\min\{1,|c|\}$, where $c$ is a suitable constant
to be fixed as $c=-2c_{0}$, with $c_{0}=\sqrt{\a}$.

The equation of motion can be written in Fourier space as
$$ i\oo\cdot\nn \left( 1 + i \e \oo\cdot\nn \right) x_{\nn} +
\e \sum_{\nn_{1}+\nn_{2}=\nn} x_{\nn_{1}} x_{\nn_{2}} =
\e f_{\nn} ,
\Eq(6.2) $$
and the formal expansion for a quasi-periodic solution
with frequency vector $\oo$ reads as
$$ x(t) = \sum_{k=0}^{\io} \e^{k} x^{(k)}(t) =
\sum_{k=0}^{\io} \e^{k} \sum_{\nn\in\zzz^{d}}
{\rm e}^{i \nn \cdot \oo t} x^{(k)}_{\nn} ,
\Eq(6.3) $$
and to see that the coefficients $x^{(k)}_{\nn}$ are well defined
to all orders $k\ge0$ one can proceed as in Section \secc(2),
with no extra difficulty. In particular the Diophantine condition
\equ(6.1) is sufficient to assure analyticity in $t$ of
the coefficients $x^{(k)}(t)$.

Also the graphical representation can be worked out as in
Section \secc(3). The only difference is that now the propagators
of the lines with non-vanishing momentum $\nn_{\ell}$,
which is defined according to \equ(3.2),
with the vectors replacing the scalars,
are given by $1/(i\oo\cdot\nn_{\ell})$, the
node factors associated with the vertices $v$ with $s_{v}=1$
are given by $F_{v}=-(i\oo\cdot\nn_{\ell_{v}})^{2}$,
and the node factors associated to the black bullets $v$ are given by
$F_{v}=f_{\nn_{v}}$, with $\nn_{v}\in\ZZZ^{d}\setminus\{\V0\}$.
All the other notations remain unchanged.

This yields that the propagators and the node factors can be bounded
as in \equ(4.2) and \equ(4.3), with just a few differences
of notation. More precisely one has
$$ \eqalign{
& \Big| \prod_{v \in V_{1}(\th) } F_{v} \Big| \,
\Big| \prod_{\ell \in L_{1}(\th) } g_{\ell} \Big|
\le \prod_{\ell \in L_{1}(\th) } \left| \oo \right|
\left| \nn_{\ell} \right| , \qquad
\Big| \prod_{\ell \in L_{1}(\th) } g_{\ell} \Big| \le
\prod_{\ell \in L_{1}(\th) } {1 \over | \oo\cdot\nn_{\ell} |}
\le C_{0}^{-1} |\nn_{\ell}|^{\t} , \cr
& \Big| \prod_{v \in V_{0}(\th) } F_{v} \Big|
\le \left( {1 \over 2c_{0}} \right)^{|V_{0}(\th)|} , \qquad
\Big| \prod_{v \in E_{W}(\th) } F_{v} \Big|
\le c_{0}^{|E_{W}(\th)|} , \cr
& \Big| \prod_{v \in E_{B}(\th) } F_{v} \Big|
\le F^{|E_{B}(\th)|}
\prod_{v \in E_{B}(\th) } {\rm e}^{-\x|\nn_{v}|} , \cr}
\Eq(6.4) $$
where the only bound which introduces a real difficulty
with respect to the case of periodic forcing terms is
the second one in the first line. Indeed it is
the source of a small divisors problem, which can not be
set only through the Diophantine condition \equ(6.1).

To each order $k$ we obtain for $x^{(k)}(t)$ a bound like
$A B^{k} k!^{\max\{1,\t\}}$, where the factor $1$ arises
from the propagators of the lines in $L_{1}(\th)$ and
the factor $\t$ from those of the lines in $L_{2}(\th)$ in \equ(6.3).
The last assertion is easily proved by reasoning as in \equ(4.3),
with $\max\{|\oo|\,|\nn_{\ell}|,C_{0}^{-1}|\nn_{\ell}|^{\t}\}\le \max
\{C_{0}^{-1},|\oo|\}|\nn_{\ell}|^{\max\{1,\t\}}$ replacing $\n_{\ell}$.
In particular only for $d=2$ and $\t=1$ we obtain the same
bound proportional to $k!$ as in the case of periodic solution
(of course with different constants $A$ and $B$).
Note that the vectors satisfying the Diophantine condition \equ(6.1)
with $\t=1$ for $d=2$ is of zero measure but everywhere dense.
An example of vector of this kind is $\oo=(1,\g_{0})$,
where $\g_{0}= (\sqrt{5}-1)/2$ is the golden section.

However, to deal with the problem of accumulation of small
divisors and discuss the issue of convergence of the series,
we need Renormalization Group techniques. The first step
is just to introduce a multiscale decomposition of the
propagators, and this leads naturally to the introduction
of clusters and self-energy graphs into the trees.
The discussion can be performed either as in
Ref.~\cita{GG1} or as in Refs.~\cita{G} (and in Ref.~\cita{GG2}).
We choose to follow Ref.~\cita{G}, which is more similar
to the present problem because the propagators
are scalar quantities and not matrices. In any case,
with respect to the quoted reference, we shall use
a multiscale decomposition involving only the quantities
$|\oo\cdot\nn_{\ell}|$, that is without introducing any dependence on $\e$
in the compact support functions. Indeed this is more suitable
to investigate the analyticity properties in $\e$,
and, as we shall see, we shall not need to exclude any
real value of $\e$ in order to give a meaning to the
resummed series, a situation more reminiscent of
Ref.~\cita{GG1} than of Refs.~\cita{G}.

In the following we confine ourselves to outline the main
differences with respect to Ref.~\cita{G}.
Let us introduce the functions $\psi_{n}$ and $\chi_{n}$,
for $n\ge0$, as in Ref.~\cita{G}, Section 5. In particular
$\psi_{n}(|x|) \neq 0$ implies $|x| \ge 2^{-(n+1)}C_{0}$ and
$\chi_{n}(|x|) \neq 0$ implies $|x| \le 2^{-n}C_{0}$.
We shall define recursively the {\it renormalized propagators}
$g^{[n]}_{\ell}=g^{[n]}(\oo\cdot\nn_{\ell};\e)$ and the {\it counterterms}
$\MM^{[n]}(\oo\cdot\nn;\e)$ on scales $n$ as
$$ \eqalign{
g^{[-1]}(x;\e) & = 1 , \qquad M^{[-1]}(x;\e) = 0 , \cr
g^{[0]}(x;\e) & = {\psi_{0}(|x|) \over ix (1+i\e x)} , \qquad
M^{[0]}(x;\e) = \sum_{k=1}^{\io} \sum_{T \in \SS^{\RR}_{k,0}}
\VV_{T}(x;\e) , \cr
g^{[n]}(x;\e) & = {\chi_{0}(|x|) \ldots \chi_{n-1}(|x|)\psi_{n}(|x|)
\over ix (1+i\e x) + \MM^{[n-1]}(x;\e)} , \cr
\MM^{[n]}(x;\e) & = \MM^{[n-1]}(x;\e) +
\chi_{0}(|x|) \ldots \chi_{n-1}(|x|)\chi_{n}(|x|) M^{[n]}(x;\e) , \cr
M^{[n]}(x;\e) & = \sum_{k=1}^{\io} \sum_{T \in \SS^{\RR}_{k,n}}
\VV_{T}(x;\e) , \cr}
\Eq(6.5) $$
where the set of renormalized self-energy graphs $\SS^{\RR}_{k,n}$
and the self-energy graphs $\VV_{T}(x;\e)$ are defined as in
Ref.~\cita{G}, Section 6. We have explicitly used that
the first contribution to the self-energy graphs is of
order $k=1$ (see Figure \graf(6.1)). Note that one
has $\chi_{0}(|x|) \ldots \chi_{n-1}(|x|)\chi_{n}(|x|)=\chi_{n}(|x|)$,
so that if $g^{[n]}(x;\e) \neq 0$ then
one has $2^{-(n+1)}C_{0} \le |x| \le 2^{-(n-1)}C_{0}$.

\midinsert
\*
\insertplotttt{226pt}{62pt}{
\ins{71pt}{58pt}{$(0)$}
\ins{89pt}{23pt}{$+$}
\ins{198pt}{42pt}{$\nn$}
\ins{190pt}{23pt}{$-\nn$}
\ins{216pt}{74pt}{$(1)$}
\ins{216pt}{35pt}{$(1)$}
}
{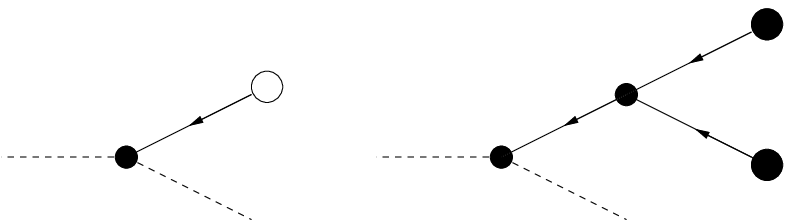}
\*
\line{\vtop{\line{\hskip1.5truecm\vbox{\advance\hsize by -3.1 truecm
\0{\nota
\figura(6.1).
First orders contributions to the counterterm arising from
self-energy graphs of order $k=1$ and $k=3$. The (dashed)
external lines do not enter into the definition of self-energy graph,
and they have been drawn only with the aim of helping visualizing
the structure of the self-energy graph.}
} \hfill} }}\*
\endinsert

Then one defines for $k\ge 1$
$$ x^{[k]}_{\nn} = \sum_{\th \in \TT_{k,\nn}} \Val(\th) , \qquad
x^{[k]}_{\V0} \= c_{k} =
{\mathop{\sum}_{\th\in\TT_{k,\V0}}}^{*} \Val(\th) ,
\Eq(6.6) $$
where the tree value is defined as
$$ \Val(\th) = \left(\prod_{\ell \in L(\th)} g^{[n_{\ell}]}_{\ell}
\right) \left(\prod_{v \in E(\th) \cup V(\th)} F_{v} \right) ,
\Eq(6.7) $$
and, as before, $*$ means that there are two lines entering the last
vertex $v_{0}$ of $\th_{0}$ and neither one exits from
an endpoint $v$ with order label $k_{v}=0$. For $k=1$
the second of \equ(6.7) has to be interpreted as $c_{1}=0$.

Furthermore one has
$$ \eqalign{
M^{[0]}(x;\e) & = M^{[0]}(0;\e) + O(\e^{2}x) , \cr
M^{[0]}(0;\e) & = - 2 \e c_{0} + M^{[0]}_{2} (0;\e) , \qquad
M^{[0]}_{2} (0;\e) = O(\e^{2}) , \cr}
\Eq(6.8) $$
and an easy computation shows (cf. Figure \graf(6.1)) that
$$ M^{[0]}_{2} (0;\e) = \e^{3} {1 \over c_{0}} \sum_{\nn\neq\V0}
\psi_{0}^{2}(|\oo\cdot\nn|) { |f_{\nn}|^{2} \over
(\oo\cdot\nn)^{2} (1 + (\e\oo\cdot\nn )^{2})} + O(\e^{4}) ,
\Eq(6.9) $$
so that in fact one has $M^{[0]}_{2} (0;\e) = O(\e^{3})$.

Moreover to higher scales one has $M^{[n]}(x;\e) =
M^{[n]}(0;\e) + O(\e^{3}x)$, with
$$ M^{[n]}(0;\e) = - \e^{3} {1 \over c_{0}} \sum_{\nn\neq\V0}
\sum_{n_{1}+n_{2}=n} \psi_{n_{1}}(|\oo\cdot\nn|)
\psi_{n_{2}}(|\oo\cdot\nn|) { |f_{\nn}|^{2} \over
(\oo\cdot\nn)^{2} (1 + (\oo\cdot\nn)^{2}) } + O(\e^{4}) ,
\Eq(6.10) $$
so that each $M^{[n]}(0;\e)$ is a higher order correction to
$M^{[0]}(0;\e)$ and it decays exponentially in $n$
(because of the compact support functions).

The following result holds.

\*

\0{\bf Lemma 5.} {\it Assume that the renormalized propagators up
to scale $n-1$ can be bounded as
$$ \left| g^{[n_{\ell}]}_{\ell} \right| \le
C_{1}^{-1} 2^{\b n_{\ell}}
\Eq(6.11) $$
for some positive constants $C_{1}$ and $\b$.
Then for all $n' \le n-1$ the number $N_{n'}(\th)$
of lines on scale $n'$ in $\th$ is bounded by
$$ N_{n'}(\th) \le K 2^{-n'/\t} \sum_{v\in E_{B}(\th)} |\nn_{v}| ,
\Eq(6.12) $$
for some positive constant $K$. If $|\e| < \e_{0}$,
with $\e_{0}$ small enough, then for all $n' \le n$ one has
$$ | M^{[n']}(x;\e) | \le D_{1} |\e|^{3} e^{- D_{2} 2^{n'/\t}} , \qquad
| \dpr_{x} M^{[n']}(x;\e) | \le
D_{1} |\e|^{3} e^{- D_{2} 2^{n'/\t}} ,
\Eq(6.13) $$
for some $C_{1}$-independent positive constants $D_{1}$ and $D_{2}$.
Only the constant $D_{1}$ depends on $\b$. The constant $\e_{0}$
can be written  as $\e_{0}=C_{1}C_{2}^{-\b}C_{3}$, with $C_{2}$ and
$C_{3}$ two positive constant independent of $\b$ and $C_{1}$.}

\*

\0{\it Proof.} The proof can be easily adapted from the proofs
of Lemma 1 and Lemma 2 of Ref.~\cita{G}. \qed

\*

So we are left with the problem of proving that the renormalized
propagators satisfy the bounds \equ(6.11). To this end
let us introduce the notation
$$ F(x) = F_{0}(x) + c_{1}(\e)\,\e + c_{2}(\e,x)\e^{2}x , \qquad
F_{0}(x) = ix (1+i\e x) ,
\Eq(6.14) $$
with $x=\oo\cdot\nn$ and the functions $c_{1}(\e)$
and $c_{2}(\e,x)$ such that $c_{1}(\e) = c + c_{3}(\e)\,\e$,
with $c\neq 0$, and the functions $|c_{2}(\e,x)|$ and $|c_{3}(\e)|$
bounded by a constant $c'$ uniformly (in $\e$ and $x$).
Recall that $\g=\min\{1,|c|\}$ and $C_{0}<\g/2$.

Fix $\l\in[0,1]$. Set $\BB_{R}(0)=\{\e\in \CCC : |\e| < R\}$ and
$\DD_{R,\l} = \{ \e = a + ib \in \BB_{\l R}(0) : |a| \ge \l |b| \}$
(see Figure \graf(6.2)). The following result refines Lemma 4.

\midinsert
\*
\insertplotttt{343pt}{109pt}{
\ins{-15pt}{60pt}{(a)}
\ins{182pt}{60pt}{(b)}
\ins{134pt}{50pt}{${\rm Re}\,\e$}
\ins{40pt}{108pt}{${\rm Im}\,\e$}
\ins{334pt}{50pt}{${\rm Re}\,\e$}
\ins{239pt}{108pt}{${\rm Im}\,\e$}
}
{fishfig14}
\*
\line{\vtop{\line{\hskip1.5truecm\vbox{\advance\hsize by -3.1 truecm
\0{\nota
\figura(6.2).
Region $\DD_{R,\l}$ in the complex $\e$-plane
for $\l=\tan \p/6$ (a) and for $\l=1$ (b).
One can write $\l=\tan\f$, where $\f$ is the angle between
the imaginary axis and the line $a=\l b$.}
} \hfill} }}
\*
\endinsert

\*

\0{\bf Lemma 6.} {\it Given $0<R<1/4C_{0}$ let $\CC_{R}$ be defined as
in Lemma 4. For all $\e\in \CC_{R}$ and all $x$ one has
$|F_{0}(x)| \ge \min\{C_{0},|x|\}/2$, while for all $\e\in\DD_{R,\l}$
one has $|F_{0}(x)| \ge \l |x|/2$.}

\*

\0{\it Proof.} Write $\e=a+ib$, so that $|F_{0}(x)|=|x|
\sqrt{(1-bx)^{2}+(ax)^{2}}$.
For $\e\in\CC_{R}$ set $A=1-\sqrt{C_{0}R}$.
If $|x|\ge C_{0}$, for $|1-bx| \le A$ one has
$|F_{0}(x)|\ge |ax^{2}| \ge b^{2}x^{2}/2R \ge C_{0}/2$, while for
$|1-bx|\ge A$ one has $|F_{0}(x)|\ge A|x| \ge |x|/2 \ge C_{0}/2$.
If $|x|\le C_{0}$, for $|1-bx| \le A$
one has $|F_{0}(x)|\ge |ax^{2}| \ge C_{0}/2 \ge |x|/2$,
while for $|1-bx|\ge A$ one has $|F_{0}(x)|\ge A|x| \ge |x|/2$. 
For $\e\in\DD_{R,\l}$ set $A=1/2$:
one finds $|F_{0}(x)|\ge \l|x|/2$. \qed

\*

Then the following result holds.

\*

\0{\bf Lemma 7.} {\it Set $x=\oo\cdot\nn$ and assume $|x|\le C_{0}$.
Then if $R$ is small enough one has
$|F(x)| \ge \l \g |x|/8$ for all $\e\in \DD_{R,\l}$.}

\*

\0{\it Proof.} Set $F_{1}(x)=F_{0}(x)+c\e$ and $\e=a+ib$.
Then $F_{1}(x)=i(x+b(c-x^{2}))+a(c-x^{2})$, and
$|F(x)| \ge |F_{1}(x)| - c'|\e|^{2}(1+|x|)$.
If $|x+b(c-x^{2})|\ge |x|/2$ and $|bc| \ge 4|x|$ one has
$|F_{1}(x)|\ge |c|\sqrt{b^{2}+a^{2}}/2\=|c\e|/2$,
so that $|F(x)|\ge |c\e|/4 \ge |cb|/4 \ge |x|$.
If $|x+b(c-x^{2})|\ge |x|/2$ and $|bc| \le 4|x|$ one has
$|F_{1}(x)|\ge \g\max\{\sqrt{x^{2}+a^{2}},|\e|/4\}/2$,
so that $|F(x)|\ge \g \sqrt{x^{2}+a^{2}}/4 \ge \g |x|/4$.
If $|x+b(c-x^{2})|\le |x|/2$ one has
$|b(c-x^{2})| \ge |x|/2$ and $|bc| \le 3 |x|$, which give
$|\e|^{2} \le 3|\e|\sqrt{a^2+x^2}/\g \le 3 \l R(|a|+|x|)/\g$,
and $|F_{1}(x)| \ge |a(c-x^{2})| \ge
|a(c-x^{2})|/2+(\l|x|/2)/2 \ge \g \l (|a|+|x|)/4$,
so that $|F(x)| \ge \g \l |x|/8$. \qed

\*

Then we can come back to the bounds of the renormalized
propagators, and prove the following result.

\*

\0{\bf Lemma 8.} {\it If $R$ is small enough
for all $n\ge 0$ and all $\e\in \DD_{R,\l}$
the renormalized propagators $g^{[n]}(x;\e)$ satisfy the bounds
\equ(6.11) with $\b=1$ and $C_{1}=\l C_{4}$,
with a $\l$-independent constant $C_{4}$.}

\*

\0{\it Proof.} The proof can be done by induction on $n$.
For $n=0$ the bound is trivially satisfied by Lemma 6.
Assuming that the bounds hold for all $n'<n$ then
we can apply Lemma 5 and deduce the bounds \equ(6.13).
In turns this implies that the renormalized propagators
on scale $n$ can be written as $g^{[n]}(x;\e)=1/F(x)$, with $F(x)$
written as in \equ(6.14) for $c=-2c_{0}$ (cf. \equ(6.8)),
and for suitable functions $c_{1}(\e)$ and $c_{2}(\e,x)$, depending
on $n$ and satisfying the properties listed after \equ(6.14)
for some $n$-independent constant $c'$. Then by Lemma 7
the renormalized propagators $g^{[n]}(x;\e)$ satisfy the same
bounds \equ(6.11) with $C_{1}=O(\l)$ for $\e\in \DD_{R,\l}$. \qed

\*

Of course for real $\e$ the bound \equ(6.11) is trivially
satisfied, with $C_{1}=2^{-1}C_{0}$. This follows
from Lemma 8 with $\l=1$, but it is obvious independent
of that result because one has $c_{1}(\e)=c\e + O(\e^{2})$,
with $c=-2c_{0} \in\RRR$. If we want to take also complex values of $\e$,
we have analyticty in a domain $D$ which can be written as
$\DD_{R} = \cup_{\l\in[0,1]} \DD_{R,\l}$.
One can easily realize that the region $\CC_{R}$ is contained
inside the domain $\DD_{R}$ (cf. Figure \graf(6.3)).
Fix $\l=\tan\f$, with $\f\in[0,\pi/4]$ (see Figure \graf(6.2)):
for all such $\f$ the line which forms an angle $\f$ with the
imaginary axis (see Figure \graf(6.2)) and passes through the
origin intersects the boundary of $\DD_{R}$ at a distance $R\tan\f$
from the origin and the boundary of $\CC_{R}$ at a distance $R\sin\f$.
Hence we have an analyticity domain of the same form as in
the case of periodic forcing terms. Nevertheless the results
found so far do not allow us to obtain Borel summability,
notwithstanding a circular analyticity domain $\CC_{R}$ is found,
as the bounds which are satisfied inside the region $\CC_{R}$
are not uniform in $\e$ (because of the dependence on $\l$).

\midinsert
\*
\insertplotttt{181pt}{109pt}{
\ins{171pt}{50pt}{${\rm Re}\,\e$}
\ins{66pt}{108pt}{${\rm Im}\,\e$}
}
{fishfig15}
\*
\line{\vtop{\line{\hskip1.5truecm\vbox{\advance\hsize by -3.1 truecm
\0{\nota
\figura(6.3).
Regions $\DD_{R}$ and $\CC_{R}$ in the complex $\e$-plane.
The first one is the grey region,
while the second one is the region contained
inside the two circles.}
} \hfill} }}
\*
\endinsert

Note hat $\b=1$ in \equ(6.11) is the same exponent appearing in the
bounds of the propagators in the formal expansion.
To obtain uniform bounds in a domain $\CC_{R}$, for some value of $R$,
we have to allow larger values of $\b$.
The following result is obtained.

\*

\0{\bf Lemma 9.} {\it Set $x=\oo\cdot\nn$ and assume $|x|<C_{0}$.
If $R$ is small enough one has $|F(x)|>\g|x|^{2}/2$ for
all $\e\in \CC_{R}$.}

\*

\0{\it Proof.} Set $F_{1}(x)=F_{0}(x)+c\e$ and $\e=a+ib$.
If $|x+b(c-x^{2})|\le |x|/2$ one has
$|b(c-x^{2})| \ge |x|/2$ and $3|x| \ge |bc| \ge |x|/4$.
Hence $|a(c-x^{2})| \ge b^{2}|c-x^{2}|/2R \ge |x|^{2}/16R|c|$,
so that one has $|F_{1}(x)| \ge |a(c-x^{2})| \ge
|ac|/4+|a(c-x^{2})|/2 \ge \g (|a|+x^{2}/16Rc)/2$.
On the other hand one has $|\e|^{2} = a^{2} + b^{2} \le
a^{2} + 9x^{2}/c^{2}$, so that $|F(x)| \ge |F_{1}(x)|-
2c'|\e|^{2} \ge |F_{1}(x)|/2 \ge \g x^{2}/2$.
The case $|x+b(c-x^{2})|\ge |x|/2$ can be discussed as in Lemma 7,
and it gives $|F(x)| \ge \g |x|/4$. \qed

\*

Then we can prove the following result
by proceeding exactly as in the proof of Lemma 8.

\*

\0{\bf Lemma 10.} {\it If $R$ is small enough
for all $n\ge 0$ and all $\e\in \CC_{R}$
the renormalized propagators $g^{[n]}(x;\e)$ satisfy the bounds
\equ(6.11) with $\b=2$ and $C_{1}$ a suitable constant.}

\*

The advantage of Lemma 8 with respect to Lemma 10
is that the bound of $R$ is better, which means that the
domain $\CC_{R}$ contained inside $\DD_{R}$ in the first case 
is larger than the domain $\CC_{R}$ of the second case.
The advantage of Lemma 10 is that it allows to obtain uniform bounds
inside the corresponding domain $\CC_{R}$. Nevertheless, because of the
factor $\b=2$, a bound $AB^{k}k!^{2\t}$ is obtained for
the coefficients $x^{(k)}(t)$ of the formal solution,
and a result analogous to Proposition 2 can be proved also for
the present case, with $N!^{2\t}$ replacing $N!$;
we do not give the details as the proof is identical.
Hence the bounds that we have are not good enough to obtain
Borel-summability in the case of quasi-periodic forcing terms,
a situation strongly reminiscent of that encountered in Ref.~\cita{GG1}.
In fact at best one can set $\t=1$ for $d=2$ (which, as noted above,
corresponds to a set of Diophantine vectors of zero measure
but everywhere dense), but this in turn implies a bound proportional
to $N!^{2}$, which is not enough to apply Nevanlinna's theorem.

The conclusion is that the resummed series
$$ x(t) = \sum_{k=0}^{\io} \m^{k} x^{[k]}(t) ,
\Eq(6.15) $$
where the coefficients $x^{[k]}(t)$ are given by
$$ x^{[k]}(t) \sum_{\nn\in\zzz^{d}} {\rm e}^{i\nn\cdot\oo t}
x^{[k]}_{\nn} ,
\Eq(6.16) $$
with $x^{[k]}_{\nn}$ defined by \equ(6.6), is well defined
and converges. In general it is not obvious -- even if expected, --
that \equ(6.15) solves the equation of motion \equ(1.1).
Indeed, unlike the case of periodic forcing terms,
we have no result, such as Nevanlinna's theorem on Borel summability, 
which we can rely upon in order to link the resummed series
to the formal series. Therefore we have to check by hand
that by expanding in powers of $\e$ the resummed series
we recover the formal power series \equ(6.3). This means
that the resummed series, which in principle could be unrelated to
the equation of motion (because of the way it has been defined),
in fact solves such an equation. Such a property
can be proved by reasoning as in Ref.~\cita{G}, Section 8.
Again we omit the details, which can be easily worked out.

We can summarize our results in the following statement.

\*

\0{\bf Theorem 3.} {\it Given the equation \equ(1.1) with $f$
analytic in its argument and $\oo$ satisfying the Diophantine
condition \equ(6.1), there exists $\e_{0}$ such that for all real
$\e$ with $|\e|<\e_{0}$ there is a quasi-periodic solution
with the same frequency vector of the forcing term. Such a solution
extends to a function analytic in a domain $\DD_{R}$ like in
Figure \graf(6.3), with $R=\e_{0}$.}

\*

The conclusion is that the summation criterion described here gives
a well defined function, which is quasi-periodic and solves
the equation of motion \equ(1.1),
but the criterion is not equivalent to Borel summability any more.
In particular the issue if such quasi-periodic solutions
are unique or not remains open, as in Ref.~\cita{GG1}.

\*\*
\section(7,Extension to more general nonlinearities)

\0When considering the equation \equ(1.5) the formal analysis
of Section \secc(2) (and of Section \secc(6) in the case
of quasi-periodic forcing terms) can be performed essentially
in the same way. If we write
$$ \eqalign{
g(x) & = \sum_{p=0}^{\io} {1 \over p!} g_{p} \,
(x-c_{0})^{p} , \qquad g_{p} =
{{\rm d}^{p} g \over {\rm d}x^{p}} (c_{0}) , \cr
[g(x)]^{(k)}_{\n} & = \sum_{p=0}^{\io} {1\over p!} g_{p}
\sum_{{\scriptstyle k_{1}+\ldots+k_{p}=k} \atop
{\scriptstyle \nn_{1}+\ldots+\nn_{p}=\nn} }
x^{(k_{1})}_{\nn_{1}} \ldots x^{(k_{p})}_{\nn_{p}} , \qquad
k \ge 0 , \cr}
\Eq(7.1) $$
then the recursive equations for $\nn\neq\V0$ are
$$ \eqalign{
x^{(0)}_{\nn} & = 0 , \cr
x^{(1)}_{\nn} & = {f_{\nn} \over i\oo\cdot\nn} , \cr
x^{(k)}_{\nn} & = - (i\oo\cdot\nn)\, x^{(k-1)}_{\nn} -
{1 \over i\oo\cdot\nn} [g(x)]^{(k-1)}_{\nn} ,
\qquad k \ge 2 , \cr}
\Eq(7.2) $$
while the compatibility condition becomes
$[g(x)]^{(k)}_{\V0} = f_{\V0}\d_{k,0}$ for $k\ge 0$.
The latter for $k=0$ gives $g(c_{0})=f_{0}$, while for $k\ge1$ gives
$ g'(c_{0}) \, c_{k} + R(c_{0},c_{1},\ldots,c_{k-1}) = 0$,
where the function $ R(c_{0},c_{1},\ldots,c_{k-1})$
depend on the coefficients to all orders $k'<k$,
hence, in particular, on the constants $c_{0},\ldots,c_{k-1}$.
Therefore the constants $c_{k}$ can be fixed iteratively as
$$ c_{k} = -{1\over g'(c_{0})} R(c_{0},c_{1},\ldots,c_{k-1}) ,
\Eq(7.3) $$
provided that one has $g'(c_{0})\neq0$, so that under the
conditions \equ(1.6) one has the formal solubility of the
equations of motion \equ(1.1). Note that the fisrt
condition in \equ(1.6) requires $f_{0} \in {\rm Ran}(g)$,
and if such a condition is satisfied then the condition
on the derivative is a genericity condition.
Note also that the class of functions
$g(x)$ which are not allowed depends on $f$ (more precisely
on its average $f_{0}$). For instance an explicit example of
function which does not satisfy \equ(1.6) is 
$g(x)=3x^{2}-2x^{3}$ if $f_{0}=1$.

The graphical representation differs from that of the previous
section as now the number of lines entering a vertex 
$v$
can assume any value $s_{v}\in\NNN$, and
if $v\notin V_{0}(\th)$ the corresponding node factor is
$$ F_{v} = - {\e \over s_{v}!} g_{s_{v}} ,
\Eq(7.4) $$
which is bounded proportionally to some constant $G$
to the power $s_{v}$. Since $\sum_{v\in V(\th)}
(s_{v}-1)=|E(\th)|-1 \le k-1$ (by Lemma 3) this
produces an overall constant $G^{2k}$ in the tree value.
Also the study of the convergence
of both the formal series and the resummed series
can then be performed as in the previous case,
and no further difficulty arises.
The constant $c$ appearing after \equ(6.14) becomes
$-g'(c_{0})$, instead of $-2c_{0}$, so that still one has
$c\neq0$ by the assumption \equ(1.6).

\*\*
\0{\bf Acknowledgments.} We are indebted to Ugo Bessi
for many enlightening discussions.


\rife{BG1}{1}{
A. Berretti, G. Gentile,
{\it Scaling properties for the radius of convergence
of a Lindstedt series: the standard map},
J. Math. Pures Appl. (9) {\bf 78} (1999), no. 2, 159--176. }
\rife{BG2}{2}{
A. Berretti, G. Gentile,
{\it Scaling properties for the radius of convergence
of a Lindstedt series: generalized standard maps},
J. Math. Pures Appl. (9) {\bf 79} (2000), no. 7, 691--713. }
\rife{BMS}{3}{
A. Berretti, S. Marmi, D. Sauzin,
{\it Limit at resonances of linearizations of some
complex analytic dynamical systems},
Ergodic Theory Dynam. Systems {\bf 20} (2000), no. 4, 963--990. }
\rife{CGM}{4}{E. Caliceti, V. Grecchi, M. Maioli,
{\it The distributional Borel summability and the large
coupling $\Phi^4$ lattice fields},
Comm. Math. Phys. {\bf 104} (1986), no. 1, 163--174;
{\it Erratum},
Comm. Math. Phys. {\bf 113} (1987), no. 1, 173--176. }
\rife{DMGB}{5}{
J.H.B. Deane, L. Marsh, G. Gentile, M.V. Bartuccelli,
{\it Invariant sets, and dynamics, of the varactor equation},
in preparation. }
\rife{GG1}{6}{
G. Gallavotti, G. Gentile,
{\it Hyperbolic low-dimensional invariant tori and summation
of divergent series},
Comm. Math. Phys. {\bf 227} (2002), no. 3, 421--460. }
\rife{GG2}{7}{
G. Gentile, G. Gallavotti, 
{\it Degenerate elliptic tori},
to appear on Comm. Math. Phys. }
\rife{G}{8}{
G. Gentile,
{\it Quasi-periodic solutions for two-level systems},
Comm. Math. Phys. {\bf 242} (2003), no. 1-2, 221--250. }
\rife{H}{9}{
G.H. Hardy,
{\it Divergent Series},
Clarendon Press, Oxford, 1949. }
\rife{MS}{10}{
S. Marmi, D. Sauzin,
{\it Quasianalytic monogenic solutions of a cohomological equation},
Mem. Amer. Math. Soc. {\bf 164} (2003), no. 780, vi+83 pp. }
\rife{N}{11}{
F. Nevanlinna,
{\it Zur Theorie der Asymptotischen Potenzreihen},
Ann. Acad. Sci. Fenn. Ser. A {\bf 12} (1916), 1-18. } 
\rife{S}{12}{
A.D. Sokal,
{\it An improvement of Watson's theorem on Borel summability},
J. Math. Phys. {\bf 21} (1980), no. 2, 261--263. }

\biblio

\ciao